\documentclass[11pt]{article}
\usepackage{amsfonts}
\usepackage{amsfonts}
\usepackage{amsfonts}
\usepackage{amsfonts}
\usepackage{amsfonts}
\usepackage{amscd}
\usepackage{amsmath}
\usepackage[all]{xy}
\usepackage{latexsym}
\usepackage{amssymb}
\usepackage{amsfonts}
 \usepackage{mathrsfs}
 \usepackage[all]{xy}
 \usepackage{latexsym,amsfonts,amssymb,bm}

\begin{document}
\baselineskip=0.6truecm

\font\one=cmbx10 scaled\magstep4 \font\fif=cmti10 scaled\magstep1
\font\fiv=cmti10 scaled\magstep5 \font\ooo=cmsy10 scaled\magstep3
\font\two=cmcsc10 \font\three=cmti8 \font\ss=cmsl8 \font\bb=cmbx8
\font\four=cmbx10 scaled\magstep1 \font\six=cmti10 scaled\magstep1
\font\small=cmr8

\def\f{\noindent}
\def\c{\centerline}
\def\ol{\overline}
\def\T{\text{T}}
\def\t{\text}
\def\d{{\text{d}}}
\def\om{\omega}
\def\Om{\Omega}
\def\sub{\subset}
\def\al{\alpha}
\def\dt{\delta}
\def\ep{\varepsilon}
\def\eq{\equiv}
\def\la{\lambda}
\def\lg{\langle}
\def\rg{\rangle}
\def\e{{{\text{e}}}}
\def\R{{\Bbb R}}
\def\C{{\Bbb C}}
\def\A{{\cal A}}
\def\B{{\cal B}}
\def\vp{\varphi}
\def\pt{\partial}
\def\na{\nabla}
\def\Dt{\Delta}
\def\pma{\pmatrix}
\def\epm{\endpmatrix}
\def\Ga{\Gamma}
\def\for{\forall}
\def\si{\sigma}
\def\wt{\widetilde}
\def\tha{\theta}
\def\V{{\cal V}}
\def\be{\beta}
\def\bm{\bmatrix}
\def\ebm{\endbmatrix}
\def\ga{\gamma}
\def\us{\underset}
\def\os{\overset}
\def\1{{\pmb 1}}
\def\tri{\triangle}
\def\Si{\Sigma}
\def\U{{\cal U}}

\title{Tilting objects in the stable category of vector bundles on the weighted projective line of type
$(2,2,2,2;\lambda)$
\thanks{ Supported in part by the National Natural
 Science Foundation of China (Grant No. 11201386), the Natural Science Foundation of Fujian Province of China (Grant No.
 2012J05009) and the Fundamental Research Funds for the Central Universities of China(Grant No.
 2012121004)}}
\author{\small\footnotesize  Jianmin Chen, Yanan Lin,  Shiquan
Ruan
\thanks{  Corresponding author, $E-mail:$ shiquanruan@stu.xmu.edu.cn}\\
  {\small\footnotesize The School of Mathematical
Sciences, Xiamen
 University, Xiamen 361005,  P.R. China}
 }

\date{}
\maketitle \par \noindent {\bf Abstract:} \ We construct a tilting
object for the stable category of vector bundles on a weighted
projective line $\mathbb{X}$ of type $(2,2,2,2;\lambda)$, consisting of five rank
two bundles and one rank three bundle, whose endomorphism algebra is
a canonical algebra associated with $\mathbb{X}$ of type (2,2,2,2).

\section{Introduction}

\hspace{4mm}The notion of weighted projective lines
was introduced by Geigle and Lenzing [1] to give a geometric
treatment to canonical algebras which
were studied by Ringel [9]. A weighted projective line can be
interpreted as obtained from the usual projective line by inserting finitely
many weights. In [1], Geigle and Lenzing proved that the category of
coherent sheaves on weighted projective lines is derived equivalent
to the category of finite dimensional modules over some canonical
algebra, and the category of vector bundles on a weighted projective
line, as additive category, is equivalent to the category of graded
Cohen-Macaulay modules over its corresponding graded ring. So the
category of vector bundles on a weighted projective line carries two
"natural" exact structures in Quillen's sense [8]. In particular,
under the exact structure which is called distinguished exact,
induced from the category of graded Cohen-Macaulay modules over its
corresponding graded ring, the categroy of vector bundles is
Frobenius with the indecomposable projective-injective objects as
all  line bundles (i.e., rank 1 vector bundles)[4]. Therefore, the
stable category of vector bundles with respect to the distinguished
exact structure is a triangulated category.

Let $\mathbb{X}$ be a weighted projective line over an algebraically
closed field, $\mbox{vect}\mathbb{X}$  the category of vector
bundles on $\mathbb{X}$, and $\underline{\mbox{vect}}\mathbb{X}$ its
stable category obtained from the vector bundles by factoring out
all the line bundles. Recently, Kussin, Lenzing and Meltzer showed a
series of interesting results for stable vector bundle categories
for weighted projective lines of triple weight type[4-6].  Among
other things, they give a tilting object in
$\underline{\mbox{vect}}\mathbb{X}$ of type $(p_{1}, p_{2},p_{3})$
with endomorphism ring
$k\overrightarrow{\mathbb{A}}_{p_{1}-1}\otimes_{k}k\overrightarrow{\mathbb{A}}_{p_{2}-1}\otimes_{k}k\overrightarrow{\mathbb{A}}_{p_{3}-1}$,
where $p_{1}, p_{2},p_{3}$ are integers greater than or equal to 2,
and  $k\overrightarrow{\mathbb{A}}_{n}$ denotes the path algebra of
oriented quiver of type $\mathbb{A}_{n}$.

For the weighted projective line of type $(2,2,2,2;\lambda)$,
Kussin, Lenzing and Meltzer [4] proved, independently by K. Ueda
[11], that the triangulated categories
$\underline{\mbox{vect}}\mathbb{X}$ and
$D^{b}(\mbox{coh}\mathbb{X})$ are equivalent. Moreover,
$\underline{\mbox{vect}}\mathbb{X}$ has tilting objects. But it is
still difficult to find a tilting object in
$\underline{\mbox{vect}}\mathbb{X}$ since we can't describe an
explicit correspondence between objects in
$\underline{\mbox{vect}}\mathbb{X}$ and
$D^{b}(\mbox{coh}\mathbb{X})$. This paper constructs tilting objects
in the stable category of vector bundles on the weighted projective
line of type $(2,2,2,2;\lambda)$.

The paper is organized as follows:

In section 2, we collect basic definitions and properties. Section 3
discusses projective covers (res. injective hulls) of vector
bundles. Theorem 3.3 gives an explicit description of projective
covers (res. injective hulls) of Auslander bundles and extension
bundles. Proposition 3.4 shows the close relationships between the
middle term of an Auslander-Reiten  sequence in
$\mbox{vect}\mathbb{X}$ and the other two terms. We give a formula
to compute the slopes of vector bundles under shift in section 4 and
describe all the exceptional objects in
$\underline{\mbox{vect}}\mathbb{X}$ in section 5. Section 6 contains
the main result, Theorem 6.2. We construct a tilting object
consisting of five bundles of rank two and one bundle of rank three,
whose endomorphism algebra is a canonical algebra of type (2,2,2,2);
furthermore, we show that there doesn't exist any tilting object
consisting only of bundles of rank two such that whose endomorphism
algebra is a canonical algebra.

\section{weighted projective line of type  $(2,2,2,2;\lambda)$}

\hspace{4mm}Throughout $k$ is an algebraically closed field,
$\lambda$ is a closed point of $\mathbb{P}_{1}(k)$ different from
$0,1,\infty$, and identify equivalences with identifications.

In this section, we recall some basic definitions and properties
about weighted projective line of type  $(2,2,2,2;\lambda)$.

Let $\mathbb{L}$ be the rank 1 abelian group on generators
$\overrightarrow{x}_{1}, \overrightarrow{x}_{2},
\overrightarrow{x}_{3},\overrightarrow{x}_{4}$ with relations
$$2\overrightarrow{x}_{1}=2\overrightarrow{x}_{2}=2\overrightarrow{x}_{3}=2\overrightarrow{x}_{4}=:\overrightarrow{c}.$$
Then $\mathbb{L}$ is an ordered group whose cone of positive
elements is $\sum_{i=1}^{4}\mathbb{N}\overrightarrow{x}_{i}$, and
each $\overrightarrow{x}\in \mathbb{L}$ can be uniquely written in
normal form
$$\overrightarrow{x}=\sum\limits_{i=1}^{4}l_{i}\overrightarrow{x}_{i}+l\overrightarrow{c}, \ \
\mbox{where} \ \ 0\leq  l_{i}\leq 1 \ \ \mbox{and} \ \ l\in \mathbb{Z}.$$

In addition, if
$\overrightarrow{x}=\sum\limits_{i=1}^{4}l_{i}\overrightarrow{x}_{i}+l\overrightarrow{c}$
is in normal form, one can define $$\overrightarrow{x}\geq 0  \ \
\mbox{ if and only if}  \ \ l_{i}\geq 0 \ \ \mbox{ for} \ \
i=1,2,3,4 \ \ \mbox{ and} \ \ l\geq 0,$$ then each
$\overrightarrow{x}\in \mathbb{L}$ satisfies exactly one of the two
possibilities $\overrightarrow{x}\geq 0$ or $ \overrightarrow{x}\leq
\overrightarrow{\omega}+\overrightarrow{c},$  where
$\overrightarrow{\omega}=2\overrightarrow{c}-\sum\limits_{i=1}^{4}\overrightarrow{x}_{i}.$
The element $\overrightarrow{c}$ is called the canonical element and
$\overrightarrow{\omega}$ is called the dualizing  element of
$\mathbb{L}$.

Denote by $S$ the commutative algebra
$$S=k[X_{1},X_{2},X_{3}, X_{4}]/I=k[x_{1},x_{2},x_{3}, x_{4}],$$
where $I$ is the homogeneous ideal generated by
$f_{1}=X^{2}_{3}-(X^{2}_{2}+X^{2}_{1})$ and
$f_{2}=X^{2}_{4}-(X^{2}_{2}+\lambda X^{2}_{1})$. Then $S$ carries an
$\mathbb{L}$-grading by setting
$\mbox{deg}x_{i}=\overrightarrow{x}_{i}(i=1,2,3,4)$, i.e.,
$S=\bigoplus\limits_{\overrightarrow{x}\in
\mathbb{L}}S_{\overrightarrow{x}},$ where
$S_{\overrightarrow{x}}S_{\overrightarrow{y}}\subseteq
S_{\overrightarrow{x}+\overrightarrow{y}}$ and $S_{0}=k$.

By an $\mathbb{L}$-graded version of the Serre construction [10],
the category  $\mbox{coh}\mathbb{X}$ of coherent sheaves
on the weighted projective line $\mathbb{X}$ of type  $(2,2,2,2;\lambda)$ is given by
$$\mbox{coh}\mathbb{X}=\mbox{mod}^{\mathbb{L}}S/\mbox{mod}_{0}^{\mathbb{L}}S,$$
where $\mbox{mod}^{\mathbb{L}}S$ is the
category of finitely generated $\mathbb{L}$-graded $S$ modules and
$\mbox{mod}_{0}^{\mathbb{L}}S$ is  the category of finite dimensional
$\mathbb{L}$-graded $S$ modules.

Use the notation $\widetilde{M}$ for $M\in \mbox{mod}^{\mathbb{L}}S$
under the quotient functor, and call the process $q:
\mbox{mod}^{\mathbb{L}}S\rightarrow \mbox{coh}\mathbb{X}, M\mapsto
\widetilde{M}$ sheafification, then it is easy to see that
$\widetilde{M}(\overrightarrow{x})=\widetilde{M(\overrightarrow{x})}$.
Call  $\mathcal{O}=\widetilde{S}$ the structure sheaf of
$\mathbb{X}$, and $\mathcal{O}(\overrightarrow{x})$ a line bundle
for any $\overrightarrow{x}\in \mathbb{L}$.

{ \bf{Proposition 2.1}}([1]) \ The category $\mbox{coh}\mathbb{X}$  is a
hereditary, abelian, $k-$linear, Hom-finite, Noetherian category
with Serre duality, i.e., $\mbox{DExt}(X,Y)= \mbox{Hom}(Y, \tau X)$,
where the $k$-equivalence $\tau: \mbox{coh}\mathbb{X}\rightarrow
\mbox{coh}\mathbb{X}$ is the shift $X\mapsto
X(\overrightarrow{\omega})$. In addition,
$\mbox{coh}\mathbb{X}=\mathcal{H}_{+}\bigvee \mathcal{H}_{0}$, where
$\mathcal{H}_{+}$ denotes the full subcategory of
$\mbox{coh}\mathbb{X}$ consisting of all objects not having a simple
subobject, $\mathcal{H}_{0}$ denotes the full subcategory of
$\mbox{coh}\mathbb{X}$ consisting of all objects of finite length,
$\bigvee$ means that each indecomposable object of
$\mbox{coh}\mathbb{X}$ is either in $\mathcal{H}_{+}$ or in
$\mathcal{H}_{0}$, and there are no non-zero morphisms from
$\mathcal{H}_{0}$ to $\mathcal{H}_{+}$.

Objects in $\mathcal{H}_{+}$ are called vector bundles, and
$\mathcal{H}_{+}$ is also denoted by $\mbox{vect}\mathbb{X}$. In
particular, all objects with the form
$\mathcal{O}(\overrightarrow{x})$ for $\overrightarrow{x}\in
\mathbb{L}$ belong to $\mbox{vect}\mathbb{X}$.

More details about $\mathcal{H}_{0}$ as follows.

{ \bf{Proposition 2.2}}([1]) \ (1) \ The simple objects of
$\mbox{coh}\mathbb{X}$ are parametrized by the projective line
$\mathbb{P}_{1}(k)$, where to each point $a\neq \{\infty, 0,
1,\lambda\}$ there is associated a tube with a unique simple
$\mathcal{S}_{a}$, called ordinary simple, where to each  $a\in
\{\infty, 0, 1,\lambda\}$ there is associated a tube with two simple
objects, called exceptional simples. Moreover, each ordinary simple
$\mathcal{S}$ satisfies $\mbox{Ext}^{1}(\mathcal{S},\mathcal{S})=k$,
while each exceptional simple $S$ satisfies
$\mbox{Ext}^{1}(\mathcal{S},\mathcal{S})=0.$

(2) \ $\mathcal{H}_{0}$ is exact, abelian, uniserial category and
decomposes into a coproduct $\coprod_{a\in
\mathbb{P}_{1}(k)}\mathcal{U}_{a}$ of connected uniserial
subcategories, whose associated quivers are tubes with mouth simple
objects of $\mbox{coh}\mathbb{X}$.

As for the structure of  $\mbox{vect}\mathbb{X}$, there have

{ \bf{Proposition 2.3}}([1])  (1) \ For each
$\overrightarrow{x},\overrightarrow{y}\in \mathbb{L}$, there has
$$\mbox{Hom}(\mathcal{O}(\overrightarrow{x}),\mathcal{O}(\overrightarrow{y}))=S_{\overrightarrow{y}-\overrightarrow{x}}.$$
In particular,
$\mbox{Hom}(\mathcal{O}(\overrightarrow{x}),\mathcal{O}(\overrightarrow{y}))\neq
0$ if and only if $\overrightarrow{x}\leq\overrightarrow{y}.$

(2) \ Each vector bundle $E$ has a filtration by line bundles
$$0=E_{0}\subset E_{1}\subset E_{2}\subset\cdots\subset E_{r}=E,$$
where each factor $L_{i}=E_{i}/E_{i-1}$ is a line bundle.

Denoted by $K_{0}(\mathbb{X})$ the Grothendieck  group of
$\mbox{coh}\mathbb{X}$, then the classes
$[\mathcal{O}(\overrightarrow{x})]$ for $0
\leq\overrightarrow{x}\leq \overrightarrow{c}$ form a
$\mathbb{Z}$-basis of $K_{0}(\mathbb{X})$, and there is a
$\mathbb{Z}$-bilinear form $\langle-,-\rangle:
K_{0}(\mathbb{X})\times K_{0}(\mathbb{X})\rightarrow \mathbb{Z}$  on
$K_{0}(\mathbb{X})$ induced by
$$\langle [X],[Y]\rangle=\mbox{dim}\ \mbox{Hom}(X,
Y)-\mbox{dim}\ \mbox{Ext}^{1}(X,Y) \ \  \mbox{for objects} \ \
X,Y\in \mbox{coh}\mathbb{X},$$ which is called Euler form.

The following are two additive functions on $K_{0}(\mathbb{X})$
called rank and degree. The rank
$\mbox{rk}:K_{0}(\mathbb{X})\rightarrow \mathbb{Z}$ is characterized
by $\mbox{rk}(\mathcal{O}(\overrightarrow{x}))=1$ for
$\overrightarrow{x}\in \mathbb{L}$ and $\mbox{rk}(\mathcal{S})=0$
for each simple object $\mathcal{S}$.

The degree  $\mbox{deg}:K_{0}(\mathbb{X})\rightarrow \mathbb{Z}$ is
uniquely determined by the following properties:

(1) \
$\mbox{deg}(\mathcal{O}(\overrightarrow{x}))=\delta(\overrightarrow{x})$
for $\overrightarrow{x}\in \mathbb{L}$, where $\delta:
\mathbb{L}\longrightarrow \mathbb{Z}$ is the group homomorphism
defined on generators by $\delta(\overrightarrow{x}_{i})=1 \
(i=1,2,3,4)$;

(2) \ $\mbox{deg}(\mathcal{O})=0$;

(3) \ $\mbox{deg}(\mathcal{S})=2$ for each ordinary simple object,
and $\mbox{deg}(\mathcal{S})=1$ for each exceptional simple object.

For each $F\in \mbox{coh}\mathbb{X}$, define the slope of $F$ as
$\mu(F)=\mbox{deg}(F)/\mbox{rk}(F).$ It is an element in $\mathbb{Q}
\cup \{\infty\}$. And according to [1], each object in
$\mbox{vect}\mathbb{X}$ has rank $>0$, then the slope belongs to
$\mathbb{Q}$; and each object in $\mathcal{H}_{0}$ has rank $0$,
then the slope $\infty$.

In particular,

{ \bf{Theorem 2.4}}({\bf Riemann-Roch formula})([7])  \ For each $X,
Y\in \mbox{coh}\mathbb{X}$, there has $$\langle[X]\oplus [\tau
X],[Y]\rangle=\mbox{rk}(X)\mbox{deg}(Y)-\mbox{deg}(X)\mbox{rk}(Y)=\mbox{rk}(X)\mbox{rk}(Y)(\mu(Y)-\mu(X)).$$

Denote by $\mathcal{H}^{q}=\mbox{add}(\mbox{ind}
\mbox{coh}^{q}(\mathbb{X}))$, where $\mbox{ind}\mbox{coh}^{q}(\mathbb{X})$
consisting of all indecomposable vector bundles of slope $q
(q\in\mathbb{Q} \cup \{\infty\})$.

{ \bf{Proposition 2.5}}([7]) (1) \ $\mathcal{H}^{q}$ is an exact Abelian subcategory of
$\mbox{coh}\mathbb{X}$ which is closed under extension;

(2) \ $\mbox{coh}\mathbb{X}=\bigvee_{q\in\mathbb{Q} \cup
\{\infty\}}\mathcal{H}^{q}$;

(3) \ There are non-zero morphisms from $\mathcal{H}^{q}$ to
$\mathcal{H}^{r}$ if and only if $q\leq r$;

(4) \ $\mathcal{H}^{\infty}$ is just $\mathcal{H}_{0}$, and each
$\mathcal{H}^{q}$ is equivalent to $\mathcal{H}_{0}$.

We are interested in $\mbox{vect}\mathbb{X}$.

{\bf Theorem 2.6}([1]) Sheafification $q:
\mbox{mod}^{\mathbb{L}}S\rightarrow \mbox{coh}\mathbb{X}$ induces an
equivalence $$\mbox{vect}\mathbb{X}= \mbox{CM}^{\mathbb{L}}S, $$
where $\mbox{CM}^{\mathbb{L}}S$ consists of all $M\in
\mbox{mod}^{\mathbb{L}}S$ satisfying
$\mbox{Hom}(E,M)=0=\mbox{Ext}^{1}(E,M)$, for each simple
$\mathbb{L}-$grade $S-$module $E$.

{\bf Remark:} The category $\mbox{vect}\mathbb{X}$ is fully embedded
as an extension-closed subcategory into two different abelian
categories
$$\mbox{coh}\mathbb{X}\hookleftarrow\mbox{vect}\mathbb{X}= \mbox{CM}^{\mathbb{L}}S\hookrightarrow \mbox{mod}^{\mathbb{L}}S.$$
So $\mbox{vect}\mathbb{X}$ carries two different "natural" exact
structures.

Now we have the exact structure on $\mbox{vect}\mathbb{X}$ induced
from $\mbox{CM}^{\mathbb{L}}S$.  A sequence $\eta: 0\rightarrow X^{\prime}\rightarrow
X\rightarrow X^{\prime\prime}\rightarrow 0$ in $\mbox{vect}\mathbb{X}$ is
called distinguished exact if and only if $\mbox{Hom}(L, \eta)$ is
exact for each line bundle $L$.

By Serre duality, a sequence $\eta: 0\rightarrow X^{\prime}\rightarrow
X\rightarrow X^{\prime\prime}\rightarrow 0$ is distinguished exact if and only
if $\mbox{Hom}( \eta, L)$ is exact for each line bundle $L$.
Moreover, each distinguished exact sequence is exact in
$\mbox{coh}\mathbb{X}$.

{\bf Proposition 2.7}([6]) The category $\mbox{vect}\mathbb{X}$ with the
structure of distinguished exact is a Frobenius category, whose
indecomposable projective-injective objects are just all the line
bundles.

Due to Happel [3], the stable category of $\mbox{vect}\mathbb{X}$
with respect to the distinguished exact structure is a triangulated
category, denoted by $\underline{\mbox{vect}}\mathbb{X}$.

{\bf Proposition 2.8}([6]) (1) \ Let $X$ be a vector bundle without direct
summand which is a line bundle,  $IX$ be the injective hull of $X$.
There exists an exact sequence $0\rightarrow X\rightarrow
IX\rightarrow X^{''}\rightarrow 0$ in the Frobenius category
$\mbox{vect}\mathbb{X}$, then $X[1]=X^{\prime\prime}$ in
$\underline{\mbox{vect}}\mathbb{X}$;

(2) \ The stable category $\underline{\mbox{vect}}\mathbb{X}$ is
triangulated, Hom-finite, Krull-Schmidt $k-$category with Serre
duality $\underline{\mbox{Hom}}(X,
Y[1])=\mbox{D}\underline{\mbox{Hom}}(Y,\tau X)$, where $\tau$ is
induced by $X\mapsto X(\overrightarrow{\omega})$. Moreover,
$\underline{\mbox{vect}}\mathbb{X}$ is homologically finite, that is
, for any $X, Y\in \underline{\mbox{vect}}\mathbb{X}$,
$\underline{\mbox{Hom}}(X,Y[n])=0,$ for $|n|\gg 0;$

(3) \ There is an action of the Picard group $\mathbb{L}$ on
$\underline{\mbox{vect}}\mathbb{X}$ by the shift, i.e., any
$\overrightarrow{x}\in \mathbb{L}$ sends $X\in
\underline{\mbox{vect}}\mathbb{X}$ to $X(\overrightarrow{x})$;

(4) \ The stable category $\underline{\mbox{vect}}\mathbb{X}$ has
Auslander-Reiten  sequences induced from the Auslander-Reiten
sequences in $\mbox{vect}\mathbb{X}$.

Furthermore,

{\bf Theorem 2.9}([4]) The stable category
$\underline{\mbox{vect}}\mathbb{X}$ and
$D^{b}(\mbox{coh}\mathbb{X})$ are equivalent as triangulated
categories.

 \section{Projective cover and injective hull}

\hspace{4mm} In order to describe   tilting objects in
$\underline{\mbox{vect}}\mathbb{X}$, we should firstly consider the
projective covers and injective hulls of vector bundles in
$\mbox{vect}\mathbb{X}$. We consider the indecomposable bundles of
rank two in $\mbox{vect}\mathbb{X}$ since all line bundles are zero
in $\underline{\mbox{vect}}\mathbb{X}$.

{\bf Proposition 3.1} For each indecomposable vector bundle $E$ of
rank two, there exists a line bundle $L$ and a non-split exact
sequence $$0 \longrightarrow L(\overrightarrow{\omega})
\longrightarrow E \longrightarrow L(\overrightarrow{x})
\longrightarrow 0 \ \text{with}\ 0\leq \overrightarrow{x} \leq
\overrightarrow{c}.$$

$\bf{Proof}$: Recalling the first step in the proof of Theorem 4.8
in [6], and noticing that $2\overrightarrow{\omega}=0$ for the
weighted projective line of $(2,2,2,2;\lambda)$, the result follows
easily. $\hfill\square$

Now we extend the notions of Auslander bundles and extension bundles
in the weighted projective lines of triple weight type in [6] to the
type $(2,2,2,2;\lambda)$.

{\bf Definition 3.2} Let $L$ be a line bundle on $\mathbb{X}$. We
call the middle term  of the Auslander-Reiten  sequence $ 0
\longrightarrow L(\overrightarrow{\omega}) \longrightarrow E
\longrightarrow L \longrightarrow 0 $ Auslander bundle associated
with $L$, and denote it by $E=E_{L}$. For $0\leq \overrightarrow{x}
< \overrightarrow{c}$, let $\eta_{\overrightarrow{x}}:
 0\longrightarrow L(\overrightarrow{\omega})
\longrightarrow E \longrightarrow L(\overrightarrow{x})
\longrightarrow 0 $ be a non-split exact sequence. The middle term
$E=E_{L}\langle \overrightarrow{x} \rangle$, which is uniquely
defined up to isomorphism, is called the extension bundle  with the
data $(L, \overrightarrow{x})$.

{\bf Remark} The definition of extension bundle is a little
different from [6] since
  $$\mbox{dim}\ \mbox{Ext}^{1}(L(\overrightarrow{x}),L(\overrightarrow{\omega}))=1\ \text{for}\ 0\leq \overrightarrow{x} <
\overrightarrow{c}$$ and $$\mbox{dim}\
\mbox{Ext}^{1}(L(\overrightarrow{c}),L(\overrightarrow{\omega}))=2$$
in the weighted projective line of $(2,2,2,2;\lambda)$.

Now we pay attention to the projective cover $PE$ of extension
bundles $E$. Noticing that $0<
\overrightarrow{x}<\overrightarrow{c}$ implies
$\overrightarrow{x}=\overrightarrow{x}_{i}$ for some $i=1,2,3,4$, we
only need to consider the cases of $E_{L}$ and $E_{L}\langle
\overrightarrow{x}_{i} \rangle, 1\leq i\leq 4$.

{\bf Theorem 3.3} Let $L$ be a line bundle,

(1) if  $E=E_{L}$, then $PE=L(\overrightarrow{\omega}) \bigoplus
(\bigoplus \limits_{i=1}^{4}L(-\overrightarrow{x}_{i}))$;

(2) if  $E=E_{L}\langle \overrightarrow{x}_{i} \rangle$, then
$PE=L(\overrightarrow{\omega}) \bigoplus (\bigoplus \limits_{j=1,
j\neq i}^{4} L(\overrightarrow{x}_{i}-\overrightarrow{x}_{j}))$.

$\bf{Proof}$: (1) Applying
$\mbox{Hom}(L(-\overrightarrow{x}_{i}),-)$ to the exact sequence
$$\xymatrix@C=0.5cm{
  0 \ar[r] & L( \overrightarrow{\omega} ) \ar[r]^{\alpha} & E \ar[r]^{\beta} & L   \ar[r] & 0 },$$
we have $\mbox{dimHom}(L(-\overrightarrow{x}_{i}),E)=\mbox{dim
Hom}(L(-\overrightarrow{x}_{i}),L)=\mbox{dim
}S_{\overrightarrow{x}_{i}}=1.$ Assuming
$\mbox{Hom}(L(-\overrightarrow{x_{i}}),E)=\langle\varphi_{i}\rangle$,
then $ \beta\varphi_{i}=x_{i} $. For each $\varphi\in
\mbox{Hom}(L(\overrightarrow{x}),E)$, we claim that there exists a
morphism $\theta \in \mbox{Hom}(L(\overrightarrow{x}),E), $ which
 factors through $\bigoplus
\limits_{i=1}^{4}L(-\overrightarrow{x}_{i})$, such that
$\beta(\varphi-\theta)=0.$ In fact, if $\beta\varphi=0$, we can
choose $\theta=0;$ if $\beta\varphi\neq 0$, then
$-\overrightarrow{x}> 0$. We write $-\overrightarrow{x}$ in normal
form
$\sum\limits_{i=1}^{4}l_{i}\overrightarrow{x}_{i}+l\overrightarrow{c},$
and discuss it in two cases:

Case 1:  $l_{i}\neq 0 $ for some $i=1,2,3,4$. Then we have
$$\mbox{dim
Hom}(L(\overrightarrow{x}),L(-\overrightarrow{x}_{i}))=\mbox{dim
Hom}(L(\overrightarrow{x}),L)=l+1.$$ So
$x_{i}:L(-\overrightarrow{x}_{i})\longrightarrow L $ induces an
isomorphism
$$\mbox{Hom}(L(\overrightarrow{x}),L(-\overrightarrow{x}_{i}))=
\mbox{Hom}(L(\overrightarrow{x}),L).$$ \noindent Hence there exists
some $\theta_{i}\in
\mbox{Hom}(L(\overrightarrow{x}),L(-\overrightarrow{x}_{i}))$ such
that $$\beta\varphi=x_{i}\theta_{i}=\beta\varphi_{i}\theta_{i},$$
then we can choose $\theta=\varphi_{i}\theta_{i}.$

Case 2: $l_{i}= 0$ for each $i=1,2,3,4, $ i.e.,
$-\overrightarrow{x}=l\overrightarrow{c}.$ Then $$\mbox{dim
Hom}(L(\overrightarrow{x}),L(-\overrightarrow{x}_{i}))=l.$$ Assume
$\mbox{Hom}(L(\overrightarrow{x}),L(-\overrightarrow{x}_{i}))=\langle\theta_{i}^{t}|
1\leq t\leq l \rangle,$ then $\{x_{i}\theta_{i}^{t}| 1\leq t\leq
l\}$ are linearly independent in the space
$\mbox{Hom}(L(\overrightarrow{x}),L)$ since $ x_{i}:
L(-\overrightarrow{x}_{i})\longrightarrow L$ is injective. Moreover,
for each $j\neq i,$ we know that $ x_{j}^{2l}: L( \overrightarrow{x
})\longrightarrow L$ can not factor through
$L(-\overrightarrow{x}_{i})$, that is, $x_{j}^{2l}\not\in \langle
x_{i}\theta_{i}^{t}| 1\leq t\leq l \rangle.$ So $\{ x_{j}^{2l};
x_{i}\theta_{i}^{t}| 1\leq t\leq l \}$ forms a basis of
$\mbox{Hom}(L(\overrightarrow{x}),L)$. Hence, there exist $k_{t}\in
k$ for  $1\leq t\leq l+1$, such that
$\beta\varphi=\sum\limits_{t=1}^{l}k_{t}x_{i}\theta_{i}^{t}+
k_{l+1}x_{j}^{2l}
=\beta(\sum\limits_{t=1}^{l}k_{t}\varphi_{i}\theta_{i}^{t}+
k_{l+1}\varphi_{j}x_{j}^{2l-1}).$ In this case, we can choose
$\theta=\sum\limits_{t=1}^{l}k_{t}\varphi_{i}\theta_{i}^{t}+
k_{l+1}\varphi_{j}x_{j}^{2l-1}.$  This finishes the proof of the
claim.

Therefore,
 there exists $\psi\in  \mbox{Hom}(L(\overrightarrow{x}),L(\overrightarrow{\omega}))$
such that $\varphi-\theta=\alpha\psi.$ So
$\varphi=\theta+\alpha\psi$, i.e., $\varphi$  factors through
$L(\overrightarrow{\omega}) \bigoplus (\bigoplus
\limits_{i=1}^{4}L(-\overrightarrow{x}_{i}))$. It's easy to see that
there are no non-zero morphisms between two different direct
summands of $L(\overrightarrow{\omega}) \bigoplus (\bigoplus
\limits_{i=1}^{4}L(-\overrightarrow{x}_{i}))$. So
$PE=L(\overrightarrow{\omega}) \bigoplus (\bigoplus
\limits_{i=1}^{4}L(-\overrightarrow{x}_{i}))$.

(2) Applying
$\mbox{Hom}(L(\overrightarrow{x}_{i}-\overrightarrow{x}_{j}),-)$ to
$\xymatrix@C=0.5cm{
  0 \ar[r] & L( \overrightarrow{\omega} ) \ar[r]^{\alpha} & E \ar[r]^{\beta} & L( \overrightarrow{x}_{i})  \ar[r] & 0 }$,
we have $\mbox{dim
Hom}(L(\overrightarrow{x}_{i}-\overrightarrow{x}_{j}),E) =\mbox{dim
Hom}(L(\overrightarrow{x}_{i}-\overrightarrow{x}_{j}),L(\overrightarrow{x}_{i}))
=1.$ Assuming
$\mbox{Hom}(L(\overrightarrow{x}_{i}-\overrightarrow{x}_{j}),E)=\langle\varphi_{j}\rangle$,
then $ \beta\varphi_{j}=x_{j} $. Using similar arguments as in (1),
we can prove that for each $\varphi\in
\mbox{Hom}(L(\overrightarrow{x}),E)$, there exists a morphism
$\theta \in \mbox{Hom}(L(\overrightarrow{x}),E), $ which  factors
through $\bigoplus \limits_{j\neq
i}L(\overrightarrow{x}_{i}-\overrightarrow{x}_{j})$, such that
$\beta(\varphi-\theta)=0.$ Therefore, there exists some morphism
$\psi\in
\mbox{Hom}(L(\overrightarrow{x}),L(\overrightarrow{\omega}))$ such
that $\varphi-\theta=\alpha\psi,$ which means
$\varphi=\theta+\alpha\psi$  factors through
$L(\overrightarrow{\omega}) \bigoplus (\bigoplus \limits_{j\neq i}
L(\overrightarrow{x}_{i}-\overrightarrow{x}_{j}))$. It's easy to see
that there are no non-zero morphisms between two different direct
summands of $L(\overrightarrow{\omega}) \bigoplus (\bigoplus
\limits_{j\neq i}
L(\overrightarrow{x}_{i}-\overrightarrow{x}_{j}))$. So
$PE=L(\overrightarrow{\omega}) \bigoplus (\bigoplus \limits_{j\neq
i} L(\overrightarrow{x}_{i}-\overrightarrow{x}_{j}))$.
$\hfill\square$

 {\bf Note:}  If $E$ is an indecomposable vector bundle of rank
two but not an extension bundle,  then $E$  fits into  a non-split
exact sequence $0 \longrightarrow L(\overrightarrow{\omega})
\longrightarrow E \longrightarrow L(\overrightarrow{c})
\longrightarrow 0 $ for some line bundle $L$ and satisfies
$E=E(\overrightarrow{\omega})$. Hence from [7] we see that $E$ is a
quasi-simple object in some homogeneous tube with slope integer. We
will give the projective cover $PE$ of $E$ in section 4.

Now we will show the relationships between projective covers of the
vector bundle in the middle term of an Auslander-Reiten sequence and
the other two terms.

{\bf Proposition 3.4} Let $0 \longrightarrow
E(\overrightarrow{\omega})\longrightarrow F \longrightarrow E
\longrightarrow  0 $ be an Auslander-Reiten sequence in
$\mbox{vect}\mathbb{X}$ with $E$ indecomposable of rank greater than
or equal to two.

(1) If $E[-1]$ is an Auslander bundle, i.e., $E[-1]=E_{L}$ for some
line bundle $L$, then $PE \bigoplus PE(\overrightarrow{\omega}) =
PF\bigoplus L(\overrightarrow{\omega}).$

(2) If else, $PE \bigoplus PE(\overrightarrow{\omega})=PF$.

{\bf{Proof:}} Since
$E(\overrightarrow{\omega})[-1]=E[-1](\overrightarrow{\omega})$,
there exists an Auslander-Reiten sequence in
$\mbox{coh}(\mathbb{X})$
$$0 \longrightarrow  E(\overrightarrow{\omega})[-1] \longrightarrow F' \longrightarrow E[-1] \longrightarrow  0 .$$
Then we obtain a commutative diagram  with distinguished exact
sequences  as follows:

 $$\xymatrix{
   &0\ar[d] & 0\ar[d] & 0\ar[d] & \\
  0   \ar[r]^{ } & E(\overrightarrow{\omega})[-1] \ar[d]_{ i_{1}} \ar[r]^{ } &
  F' \ar[d]_{ i} \ar[r]^{ } &  E[-1] \ar[d]_{ i_{2}}  \ar[r]^{ } & 0   \\
  0   \ar[r]^{ } & PE(\overrightarrow{\omega}) \ar[d]_{\pi_{1}} \ar[r]^{(1\ 0)^{t} \ \ \ \ \  } &
  PE(\overrightarrow{\omega})\bigoplus PE  \ar[d]_{\pi} \ar[r]^{ \ \  \ \ \ \ (0\ 1)} &  PE  \ar[d]_{\pi_{2}}  \ar[r]^{ } &  0  \\
  0 \ar[r]^{ } &E(\overrightarrow{\omega}) \ar[d]\ar[r]^{f} & F\ar[d] \ar[r]^{g} & E \ar[d] \ar[r]^{} & 0  \\
   &0 & 0 & 0 & }$$

We know that the natural projective morphism $\overline{\pi}:PF
\longrightarrow F$  induces a morphism $\delta:
PE(\overrightarrow{\omega})\bigoplus PE\longrightarrow PF$. So we
get a commutative diagram :
$$\xymatrix{
  0   \ar[r]^{ } &  F' \ar[d]_{\delta^{'} } \ar[r]^{i \ \ \ \ \ \ \ \ } &
  PE(\overrightarrow{\omega})\bigoplus PE \ar[d]_{\delta } \ar[r]^{\ \ \ \ \ \ \ \   \pi} &   F  \ar[d]_{id}  \ar[r]^{ } &  0  \\
  0 \ar[r]^{ } &F[- 1] \ar[r]^{ } &PF \ar[r]^{\overline{\pi}} &   F \ar[r]^{} & 0   }.$$
Then from the snake lemma in $\mbox{coh}(\mathbb{X})$, we know that
$\delta^{'}$ is  surjective and there is  an isomorphism $\rho:
\mbox{Ker}(\delta')\longrightarrow \mbox{Ker}(\delta).$

  Now we have the following exact commutative diagram  :
 $$\xymatrix{
   0   \ar[r]^{ } &  \mbox{Ker}(\delta^{'}) \ar[d]_{\theta^{'} } \ar[r]^{\rho } &
   \mbox{Ker}(\delta ) \ar[d]_{\theta } \ar[r]^{ } &   0  \ar[d]_{ }  \ar[r]^{ } &  0  \\
  0   \ar[r]^{ } &  F' \ar[d]_{\delta^{'} } \ar[r]^{i \ \ \  \ \ \ } &
  PE(\overrightarrow{\omega})\bigoplus PE \ar[d]_{\delta } \ar[r]^{ \ \ \ \ \ \ \pi} &   F  \ar[d]_{id}  \ar[r]^{ } &  0  \\
  0 \ar[r]^{ } &F[- 1] \ar[r]^{ } &PF \ar[r]^{\overline{\pi}} &   F \ar[r]^{} & 0,   }$$
and $\delta$ is split surjective implies $\theta$ is split
injective. Hence, there exists a morphism $\sigma:
PE(\overrightarrow{\omega})\bigoplus PE\longrightarrow
\mbox{Ker}(\delta )$, such that $\sigma\theta=id_{\mbox{Ker}(\delta
)}$. Let $\sigma'=\rho^{-1}\sigma i: F'\longrightarrow
\mbox{Ker}(\delta^{'}) $, then $\sigma'\theta'=\rho^{-1}\sigma
i\theta'=\rho^{-1}\sigma\theta\rho=id_{\mbox{Ker}(\delta' )}$, which
implies  $\theta'$ is split injective. So we obtain a split exact
sequence:
$$\xymatrix@C=0.3cm{
  0 \ar[rr]^{} && \mbox{Ker}(\delta' ) \ar[rr]^{\theta'} && F' \ar[rr]^{\delta'} &&F[-1]  \ar[rr]^{} &&0 }.$$
(1) If  $E[-1]$ is an Auslander bundle, say $E[-1]=E_{L}$, then
there exists an  Auslander-Reiten sequence in
$\mbox{vect}\mathbb{X}$ as follows:
$$0 \longrightarrow
E[-1](\overrightarrow{\omega})\longrightarrow
F^{\prime\prime}\bigoplus L(\overrightarrow{\omega}) \longrightarrow
E[-1] \longrightarrow  0 ,$$ where $F^{\prime\prime}$ is the unique
indecomposable  object with rank three and socle $L$. Hence
$F^{\prime}=F^{\prime\prime}\bigoplus L(\overrightarrow{\omega})$
and $\mbox{Ker}(\delta^{'})=L(\overrightarrow{\omega})$.
Then  it follows that $PE \bigoplus PE(\overrightarrow{\omega})=PF\bigoplus L(\overrightarrow{\omega});$\\
(2) If else, $F'$ has no direct summands from line bundles. But $
\mbox{Ker}(\delta' )=\mbox{Ker}(\delta )$, which is a direct sum of
line bundles. It follows that $\mbox{Ker}(\delta' )=0$. Hence, $PE
\bigoplus PE(\overrightarrow{\omega})=PF$. $\hfill\square$

By duality, we have analogous results of injective hull $IE$ for $E$
as follows:

{\bf Theorem 3.5}  For each line bundle $L$,

(1) if  $E=E_{L}$, then $IE=L \bigoplus (\bigoplus
\limits_{i=1}^{4}L(\overrightarrow{\omega}+\overrightarrow{x}_{i}))$;

(2) if  $E=E_{L}\langle \overrightarrow{x}_{i} \rangle$, then
$IE=L(\overrightarrow{x}_{i}) \bigoplus (\bigoplus \limits_{j=1,
j\neq i}^{4} L(\overrightarrow{\omega}+\overrightarrow{x}_{j}))$.

{\bf Proposition 3.6} Let $0 \longrightarrow
E(\overrightarrow{\omega})\longrightarrow F \longrightarrow E
\longrightarrow  0 $ be an Auslander-Reiten sequence in
$\mbox{vect}\mathbb{X}$ with $E$ indecomposable of rank  greater
than or equal to two.

(1) If $E[-1]$ is an Auslander bundle, i.e., $E[-1]=E_{L}$ for some
line bundle $L$, then $IE \bigoplus
IE(\overrightarrow{\omega})=IF\bigoplus L.$

(2) If else, $IE \bigoplus IE(\overrightarrow{\omega})=IF$.

\section{Slopes under shift}

\hspace{4mm} To show an object is a tilting object, it is necessary
to discuss the change of slope of a vector bundle under shift. In
this section, we will give a formula to compute the slope of a
vector bundle under shift.

{\bf Lemma 4.1} Let $E$ be an indecomposable vector bundle with $
\mbox{rk}(E) \geq 2$,
$$0 \longrightarrow E(\overrightarrow{\omega}) \longrightarrow F_{1}\bigoplus F_{2} \longrightarrow E \longrightarrow 0 $$
\noindent be an Auslander-Reiten sequence with $F_{1}$
indecomposable and  $ \mbox{rk }(F_{1}) \geq 2$, then  $\mu(E[-1])
=\mu(F_{1}[-1]) $.

$\bf{Proof}$: Since
$\underline{\mbox{Hom}}(E(\overrightarrow{\omega}), F_{1})\neq 0$,
and since the shift functor  is an equivalence in
$\underline{\mbox{vect}}\mathbb{X}$, we obtain that
$\underline{\mbox{Hom}}(E(\overrightarrow{\omega})[-1],
F_{1}[-1])\neq 0$. Hence, $\mu(E(\overrightarrow{\omega})[-1]) \leq
\mu(F_{1}[-1])$. Similarly, we can obtain  $\mu(F_{1}[-1]) \leq
\mu(E[-1])$. Moreover, $ E(\overrightarrow{\omega})[-1]=
E[-1](\overrightarrow{\omega})$  induces
$\mu(E(\overrightarrow{\omega})[-1])=\mu(E[-1])$. So $\mu(E[-1])
=\mu(F_{1}[-1]) $. $\hfill\square$

{\bf Remark:} By Lemma 4.1, we know that two indecomposable objects
$E, F\in \mbox{vect}\mathbb{X}$ with rank $\geq 2$ are in the same
tube in the Auslander-Reiten quiver of $\mbox{vect}\mathbb{X}$ if
and only if $E[-1], F[-1]$ are so.

{\bf Lemma 4.2}  For each vector bundle $F$ and
$\overrightarrow{x}\in \mathbb{L}$ with
$\delta(\overrightarrow{x})=0$,

(1) if $\mu (F) > 0,$ then $\mbox{dim Hom}(E_{\mathcal
{O}(\overrightarrow{x})},F)= \mbox{deg} (F) $;

(2) if $\mu (F) < 0,$ then $\mbox{dim Hom}(F, E_{\mathcal
{O}(\overrightarrow{x})})=-\mbox{deg} (F) $.

{\bf{Proof : }} (1) If  $\mu(F)> 0$, applying   $\mbox{Hom}(-,F)$ to
$$\eta : 0 \longrightarrow
\mathcal{O}(\overrightarrow{\omega}+\overrightarrow{x})\longrightarrow
E_{\mathcal {O}(\overrightarrow{x})}\longrightarrow
\mathcal{O}(\overrightarrow{x})\longrightarrow 0$$ and by the
Riemann-Roch formula, we obtain
$$
\begin{aligned}
\mbox{dim Hom}(E_{\mathcal {O}(\overrightarrow{x})},F)&=\mbox{ dim
Hom}(\mathcal{O}(\overrightarrow{\omega}+\overrightarrow{x})\bigoplus
\mathcal{O}(\overrightarrow{x}), F)\\&
=\langle[\mathcal{O}(\overrightarrow{\omega}+\overrightarrow{x})]\bigoplus
[\mathcal{O}(\overrightarrow{x})],[F]\rangle\\
&=\mbox{rk}(\mathcal{O}(\overrightarrow{x}))\mbox{deg}(F)-\mbox{deg}(\mathcal{O}(\overrightarrow{x}))\mbox{rk}(F)\\
&= \mbox{deg}(F).
\end{aligned}$$

(2)  If $\mu(F)< 0$, applying  $\mbox{Hom}(F,-)$ to $\eta$ we obtain
$$\mbox{dim Hom}(F, E_{\mathcal
{O}(\overrightarrow{x})})=-\mbox{deg}(F).$$

$\hfill\square$

Now we can compute the slope of $E[-1]$ where $E$ is a vector bundle
of rank $\geq 2$ and slope $0$ or $\frac{1}{2}$.

{\bf Lemma 4.3} For each indecomposable object $E\in
\mbox{vect}\mathbb{X}$ with $\mbox{rk}(E) \geq 2$,

(1) if $\mu(E) =0$, then $\mu(E[-1]) =-\frac{4}{3}$;

(2) if $\mu(E)=\frac{1}{2}$, then $\mu(E[-1]) =-\frac{1}{2}$.

$\bf{Proof: }$ (1)By Lemma 4.1, we can reduce to the case
$\mbox{rk}(E)=2.$

If $E(\overrightarrow{\omega})\neq E$, then according to the
structure of the Auslander-Reiten quiver of $\mbox{vect}\mathbb{X}$,
$E$ is an Auslander bundle, i.e., there exists a line bundle $L$
with $\mu(L)=0$ such that $E=E_{L}$. Then
$PE=L(\overrightarrow{\omega}) \bigoplus (\bigoplus
\limits_{i=1}^{4}L(-\overrightarrow{x}_{i})), $ which implies
$\mbox{deg}(PE)=-4$ and $\mbox{rk}(PE)=5.$ From the distinguished
exact sequence
$$
  0 \longrightarrow E[-1]\longrightarrow PE \longrightarrow E \longrightarrow0 ,$$
we obtain $\mbox{deg}(E[-1])=-4$ and $\mbox{rk}(E[-1])=3$. Hence
$\mu (E[-1])=-\frac{4}{3}$.

If $E(\overrightarrow{\omega}) = E$, then $E$ is a quasi-simple
object in some homogeneous tube with slope zero. For the line bundle
$L=\mathcal{O}({-\overrightarrow{x}_{1}})$, each direct summand of
$IE_{L}\langle \overrightarrow{x}_{i}\rangle$ is with slope zero.
Hence for  $ 1\leq i\leq 4$, we have $\mbox{Hom}(IE_{L}\langle
\overrightarrow{x}_{i}\rangle,E)=0$, because there are no non-zero
morphisms between different tubes with the same slope. By the
Riemann-Roch formula, we obtain that
$$\begin{aligned}
&\ \ \ \ \mbox{dim Hom}(E_{L}\langle
\overrightarrow{x}_{i}\rangle\bigoplus E_{L}\langle
\overrightarrow{x}_{i}\rangle(\overrightarrow{\omega}),E)\\&=\langle
[E_{L}\langle \overrightarrow{x}_{i}\rangle]\bigoplus [E_{L}\langle
\overrightarrow{x}_{i}\rangle(\overrightarrow{\omega})],
[E]\rangle\\& =\mbox{rk}(E_{L}\langle \overrightarrow{x}_{i}\rangle)
\mbox{deg}(E)-\mbox{deg}(E_{L}\langle \overrightarrow{x}_{i}\rangle)
\mbox{rk}(E) \\&=2.
\end{aligned}$$ Then it follows that  $$\mbox{dim
Hom}(E_{L}\langle \overrightarrow{x}_{i}\rangle,E)=\mbox{dim
Hom}(E_{L}\langle
\overrightarrow{x}_{i}\rangle(\overrightarrow{\omega}),E)=1.$$
Hence, we obtain that $\mbox{dim}\
\underline{\mbox{Hom}}(E_{L}\langle
\overrightarrow{x}_{i}\rangle,E)=1$, which implies $$\mu(E[-1])\geq
\mu (E_{L}\langle \overrightarrow{x}_{i}\rangle[-1])=-\frac{3}{2}.$$
So each direct summand of $PE$ is with slope $-1$, i.e., $PE\in
\mathcal {H}^{-1}$. Moreover, for each $\overrightarrow{x}\in
\mathbb{L}$ with $ \delta(\overrightarrow{x})=-1,$ we have
$$\mbox{dim Hom}(\mathcal{O}(\overrightarrow{x})\bigoplus
\mathcal{O}(\overrightarrow{x}+\overrightarrow{\omega}),E)=\mbox{deg}
E(-\overrightarrow{x})=2.$$ It follows that $$\mbox{dim
Hom}(\mathcal{O}(\overrightarrow{x}) ,E)= \mbox{dim
Hom}(\mathcal{O}(\overrightarrow{x}+\overrightarrow{\omega})
,E)=1.$$ Hence,
$PE=\bigoplus\limits_{\delta(\overrightarrow{x})=-1}\mathcal{O}(\overrightarrow{x}).$
So $\mbox{deg}(PE)=-8,$ and $\mbox{rk}(PE)=8.$ From the
distinguished exact sequence
$$
  0 \longrightarrow E[-1] \longrightarrow PE \longrightarrow E \longrightarrow 0 ,$$
we obtain $\mbox{deg}(E[-1])=-8$ and $\mbox{rk}(E[-1])=6$.
Therefore, $\mu(E[-1])=-\frac{4}{3}$.

(2)According to the Lemma 4.1, we can reduce to the case that $E$ is
a quasi-simple object in $\mbox{vect}\mathbb{X}$.

Indeed,  if $E(\overrightarrow{\omega}) \neq E$  then
$\mbox{rk}(E)=2$, and $E=E_{L}\langle \overrightarrow{x}_{i}\rangle$
for some  line bundle $L$ with $\mu(L)=0$ and  $1\leq i\leq 4$. So
$PE=L(\overrightarrow{\omega}) \bigoplus (\bigoplus \limits_{j=1,
j\neq i}^{4} L(-\overrightarrow{x}_{j}))$. From the distinguished
exact sequence
$$ 0 \longrightarrow E[-1]\longrightarrow PE \longrightarrow E \longrightarrow 0 ,$$
\noindent we obtain $\mbox{deg}(E[-1])=-1$ and $\mbox{rk}(E[-1])=2$,
hence $\mu(E[-1]) =-\frac{1}{2}$.

If $E(\overrightarrow{\omega})= E$, then $\mbox{deg }(E)=2\
\text{and}\ \mbox{rk }(E)=4$. By the Riemann-Roch formula and
noticing that $\mu
(E_{\mathcal{O}(-\overrightarrow{x}_{1})}[1])=\frac{1}{3}$, we
obtain
$$\mbox{dim }
\underline{\mbox{Hom}}(E_{\mathcal{O}(-\overrightarrow{x}_{1})}[1],
E)=\mbox{dim Hom}(E_{\mathcal{O}(-\overrightarrow{x}_{1})}[1],
E)=1,$$ which implies that $$\mu(E[-1])\geq \mu
(E_{\mathcal{O}(-\overrightarrow{x}_{1})})=-1.$$ Moreover,
$E(\overrightarrow{\omega})[-1] = E[-1]$ implies $\mu(E[-1])\neq
-1$. Hence $PE\in \mathcal {H}^{0}$. Furthermore, for each
$\overrightarrow{x}\in \mathbb{L}$ with $
\delta(\overrightarrow{x})=0 $, by  arguments as above, we obtain
$\mbox{dim Hom}(\mathcal{O}(\overrightarrow{x}) ,E)=1.$ So we get
$PE=\bigoplus\limits_{\delta(\overrightarrow{x})=0}\mathcal{O}(\overrightarrow{x})$.
Therefore, $\mbox{deg}(E[-1])=-2$ and $\mbox{rk}(E[-1])=4$. It
follows that $\mu(E[-1]) =-\frac{1}{2}$. $\hfill\square$

{\bf Theorem 4.4}  Let $E,F$ be two non-isomorphism indecomposable
objects in $\mbox{vect}\mathbb{X}$ with rank $\geq 2$, then
$\mu(E)=\mu(F)$ if and only if $\mu(E[-1])=\mu(F[-1]) $.

$\bf{Proof:}$ If $\mu(E)=\mu(F)=n$, then
$\mu(E[-1])=\mu(E(-n\overrightarrow{x}_{1})[-1])+\delta(n\overrightarrow{x}_{1})=n+
\mu(E(-n\overrightarrow{x}_{1})[-1])=n-\frac{4}{3} =\mu(F[-1])$.

If $\mu(E)=\mu(F)=n+\frac{q}{p}$ with $p> q$, according to Lemma
4.1, we can reduce to the following two cases:

Case 1,  $\mbox{rk}(E)=\mbox{rk}(F)=p$. Notice that  $E,
E(\overrightarrow{\omega})$, and
$E(\overrightarrow{x}_{i}-\overrightarrow{x}_{j})(1\leq i< j\leq 4)$
are all non-isomorphism vector bundles with rank $p$ and slope
$n+\frac{q}{p}$, so there exists an element $\overrightarrow{x}\in
\mathbb{L}$ satisfying $ \delta(\overrightarrow{x})=0$, such that
$F=E(\overrightarrow{x})$. Hence
$\mu(F[-1])=\mu(E(\overrightarrow{x})[-1])=\mu(E[-1])+\delta(\overrightarrow{x})=\mu(E[-1])$.

Case 2, one of them, say $F$ satisfies   $\mbox{rk}(F)=2p$ and
$F(\overrightarrow{\omega}) = F$. If $\mu(E[-1])\neq\mu(F[-1]) $,
without loss of generality, we may assume  $\mu(E[-1])< \mu(F[-1])$.
Then by the Riemann-Roch formula, $\mbox{Hom}(E[-1],F[-1])\neq 0$.
But $\underline{\mbox{Hom}}(E[-1],F[-1])=\underline{\mbox{Hom}}(E ,F
)= 0$. So there exists some line bundle $L$ such that
$\mbox{Hom}(E[-1],L)\neq 0$ and $\mbox{Hom}(L,F[-1])\neq 0$, which
implies that  $\mu(E[-1])\leq \mu(L)$ and
$\underline{\mbox{Hom}}(E_{L(\overrightarrow{\omega})},F[-1])\neq
0$. Then  $\mu(E_{L(\overrightarrow{\omega})})< \mu(F[-1])$ since
$F[-1]= F[-1](\overrightarrow{\omega})$. Hence $\mu(E[-1])\leq
\mu(L)= \mu(E_{L(\overrightarrow{\omega})})< \mu (F[-1])$, which is
a contradiction.

Now we have shown that $\mu (E)=\mu (F)$ implies $\mu (E[-1])=\mu
(F[-1]) $. Dually, we can also obtain  $\mu (E[1])=\mu (F[1]) $ from
$\mu (E)=\mu (F)$. $\hfill\square$

Now we can give a formula to compute the slope of a vector bundle
under shift.

{\bf Proposition 4.5} For each indecomposable vector bundle $E$ with
$\mbox{rk}(E)\geq 2$,
assume $\mu (E) =n+\frac{q}{p}$, where $p, q, n \in \mathbb{Z}, 0\leq\frac{q}{p} <1$ and $(p,q)=1$. \\
(1) If $0\leq\frac{q}{p} \leq\frac{1}{3}$, then $\mu (E[-1]) =n-\frac{4p-11q}{3p-8q}  $.\\
(2) If $\frac{1}{3}<  \frac{q}{p}< 1$, then $\mu (E[-1])
=n+\frac{q}{p-4q}$.

$\bf{Proof:} $ By assumption, we have $\mu
(E(-n\overrightarrow{x}_{1})[-1])=\mu
(E[-1](-n\overrightarrow{x}_{1})) =\mu (E[-1])-
n\delta(\overrightarrow{x}_{1})=\mu (E[-1])- n$, so $\mu
(E[-1])=n+\mu (E(-n\overrightarrow{x}_{1})[-1])$. Moreover, $\mu
(E(-n\overrightarrow{x}_{1}))=\mu (E)-
n\delta(\overrightarrow{x}_{1})=\mu (E)- n=\frac{q}{p}$. Hence, we
only need to show the result for $n=0$. Moreover, according to
Theorem 4.4, we can assume that $E$ is a quasi-simple object in some
homogeneous tube, that is,  $\mbox{deg }(E)=2q $ , $\mbox{rk
}(E)=2p$, and $E=E(\overrightarrow{\omega})$.

(1) Since $ -\frac{4}{3}\leq \mu (E[-1]) \leq  -1$, we have
$PE=I(E[-1])\in \mathcal {H}^{0}\vee \mathcal {H}^{-1}$. Suppose
$\mbox{deg}(E[-1])=-2d\ \text{and}\ \mbox{rk}(E[-1])=2r$. For each
$\overrightarrow{x}\in \mathbb{L}$ satisfying
$\delta(\overrightarrow{x})=0$, we have $$\mbox{dim
Hom}(\mathcal{O}(\overrightarrow{x})\bigoplus
\mathcal{O}(\overrightarrow{x}+\overrightarrow{\omega}),E)=\mbox{dim
Hom}(E_{\mathcal{O}(\overrightarrow{x})} ,E)=\mbox{deg}(E)=2q.$$ It
follows that  $$\mbox{dim Hom}(\mathcal{O}(\overrightarrow{x}) ,E)=
\mbox{dim
Hom}(\mathcal{O}(\overrightarrow{x}+\overrightarrow{\omega})
,E)=q.$$ And for each $\overrightarrow{y}\in \mathbb{L}$ satisfying
$\delta(\overrightarrow{y})=-1$, we have $$\mbox{dim Hom}(E[-1],
\mathcal{O}(\overrightarrow{y})\bigoplus
\mathcal{O}(\overrightarrow{y}+\overrightarrow{\omega}))
=-\mbox{deg}(E[-1])-\mbox{rk}(E[-1])=2d-2r,$$ which implies
$$\mbox{dim Hom}(E[-1], \mathcal{O}(\overrightarrow{y}) )= \mbox{dim
Hom}(E[-1],
\mathcal{O}(\overrightarrow{y}+\overrightarrow{\omega}))=d-r.$$
Hence,
$$PE=I(E[-1])=(\bigoplus\limits_{\delta(\overrightarrow{x})=0}\mathcal{O}(\overrightarrow{x})^{q})\bigoplus
(\bigoplus\limits_{\delta(\overrightarrow{y})=-1}\mathcal{O}(\overrightarrow{y})^{d-r}).$$
It follows that $\mbox{deg}(PE)=-8(d-r)$ and $ \mbox{rk}(PE)
=8(d-r)+8q$. On the other hand, from the exact sequence
$$
  0 \longrightarrow E[-1] \longrightarrow PE \longrightarrow E \longrightarrow 0 ,$$
we have $\mbox{deg}(PE)=\mbox{deg}(E[-1])+\mbox{deg}(E)=2q-2d$ and $
\mbox{rk}(PE)=\mbox{rk}(E[-1])+\mbox{rk}(E )=2p+2r$. Thus we obtain
$d=4p-11q, r=3p-8q$. Hence,
 $\mu (E[-1]) =-\frac{4p-11q}{3p-8q}$.

(2)Since $ -1< \mu (E[-1])  < -\frac{1}{3}$, we have $PE\in \mathcal
{H}^{0}$. For each $\overrightarrow{x}\in \mathbb{L}$ satisfying
$\delta(\overrightarrow{x})=0$, we have $$\mbox{dim
Hom}(\mathcal{O}(\overrightarrow{x})\bigoplus
\mathcal{O}(\overrightarrow{x}+\overrightarrow{\omega}),E)=\mbox{deg}(E)=2q.$$
It follows that $$\mbox{dim Hom}(\mathcal{O}(\overrightarrow{x})
,E)= \mbox{dim
Hom}(\mathcal{O}(\overrightarrow{x}+\overrightarrow{\omega})
,E)=q.$$ Hence,
$$PE=\bigoplus\limits_{\delta(\overrightarrow{x})=0}\mathcal{O}(\overrightarrow{x})^{q}.$$
So $\mbox{deg}(PE)=0$ and $ \mbox{rk}(PE) =8q$. Then from the exact
sequence
$$
  0 \longrightarrow E[-1] \longrightarrow PE\longrightarrow E \longrightarrow 0 ,$$
we get $\mbox{deg}(E[-1]) =-2q$ and $\mbox{rk}(E[-1])=8q-2p,$ which
implies that $\mu (E[-1]) =\frac{q}{p-4q}$. $\hfill\square$

Dually  we obtain the following   result which is similar to the
previous one:

{\bf Proposition 4.6} For each indecomposable vector bundle $E$ with
$\mbox{rk}(E)\geq 2$,
assume $\mu (E) =n+\frac{q}{p}$, where $p, q, n\in \mathbb{Z}, 0\leq\frac{q}{p} <1$ and $(p,q)=1$. \\
(1) If $0\leq\frac{q}{p} \leq\frac{2}{3}$, then $\mu (E[1]) =n+\frac{4p-5q}{3p-4q}  $.\\
(2) If $\frac{2}{3}<  \frac{q}{p}< 1$, then $\mu
(E[1])=n+\frac{12p-19q}{5p-8q}$.

{\bf Remark:} If $E$ is an indecomposable vector bundle of rank two
and $E=E(\overrightarrow{\omega})$, then there exists a line bundle
$L$ and a non-split exact sequence $0 \longrightarrow
L(\overrightarrow{\omega}) \longrightarrow E \longrightarrow
L(\overrightarrow{c}) \longrightarrow 0 $. By the proof of
Proposition 4.5 and 4.6, we can show that
$PE=\bigoplus\limits_{\delta(\overrightarrow{x})=0}L(\overrightarrow{x})$
and
$IE=\bigoplus\limits_{\delta(\overrightarrow{x})=2}L(\overrightarrow{x})$.

{\bf Corollary 4.7} For each indecomposable object
$E\in\mbox{vect}\mathbb{X}$ with  $\mbox{rk}(E)\geq 2$, if $\mu
(E) =n\in \mathbb{Z}$ ,  then for each $m\in  \mathbb{N}$,\\
(1)  $\mu (E[m]) =n+m+\frac{m}{2m+1}$;\\
(2) $\mu (E[-m]) =2n-\mu (E[m]) =n-m-\frac{m}{2m+1}$.

$\bf{Proof:} $ (1)We prove the result by induction. For $m=1$, $\mu
(E[1])= n +\frac{4}{3}=n+1+\frac{1}{3}$. Suppose for $m=k$, $\mu
(E[k]) =n+k+\frac{k}{2k+1}$. Then for $m=k+1$, notice
$0\leq\frac{k}{2k+1}\leq\frac{2}{3}$, so $\mu (E[k+1])=\mu (E[k][1])
=n+k+\frac{4(2k+1)-5k}{3(2k+1)-4k}=n+k+1+\frac{ k+1 }{2(k+1)+1} $.
Therefore,  for each $m\in  \mathbb{N}$, we have $\mu (E[m])
=n+m+\frac{m}{2m+1}$.

(2)Analogously. $\hfill\square$

\section{Exceptional objects}

\hspace{4mm} Direct summands of   tilting objects are exceptional
objects. This section is due to discuss exceptional objects in
$\underline{\mbox{vect}}\mathbb{X}$. Theorem 5.5 shows all
exceptional objects in $\underline{\mbox{vect}}\mathbb{X}$.

 {\bf Lemma 5.1}([6],Proposition 2.7)
 If $X, Y$ are both exceptional in $\mbox{coh}\mathbb{X}$, and $[X]=[Y]$ in  $K_{0}(\mathbb{X})$, then $X=Y$.

 {\bf Proposition 5.2} For each line bundle $L,$ the vector bundle
 $E_{L}\langle \overrightarrow{x}_{i}\rangle$ is exceptional in $\mbox{coh}\mathbb{X},$
and $E_{L}\langle \overrightarrow{x}_{i} \rangle[1]=E_{L}\langle
\overrightarrow{x}_{i} \rangle(\overrightarrow{x}_{j})$, for each
$j\neq i .$

$\bf{Proof}$: For simplification, we write $E_{L}\langle
\overrightarrow{x}_{i}\rangle$ by $E$ throughout the proof. Applying
$\mbox{Hom}(L(\overrightarrow{\omega}),-)$ to the exact
$$\eta_{i}: 0 \longrightarrow L(\overrightarrow{\omega}) \longrightarrow E \longrightarrow L(\overrightarrow{x}_{i}) \longrightarrow 0 ,$$
we obtain $$\mbox{Hom}(L(\overrightarrow{\omega}),E )=k\ \text{and}\
\mbox{Ext}^{1}(L(\overrightarrow{\omega}),E)=0.$$ Similarly,
applying $\mbox{Hom}(L(\overrightarrow{x}_{i}) ,-)$ to $\eta_{i}$ we
obtain
$$\mbox{Hom}(L(\overrightarrow{x}_{i}),E)
=\mbox{Ext}^{1}(L(\overrightarrow{x}_{i}),E)=0.$$ Then applying
$\mbox{Hom}(-, E)$ to $\eta_{i}$ we obtain a long exact sequence:
$$\begin{aligned} 0& \longrightarrow  \mbox{Hom}(L(\overrightarrow{x}_{i}),E) \longrightarrow
\mbox{Hom}(E,E) \longrightarrow
\mbox{Hom}(L(\overrightarrow{\omega}),E)\\&\longrightarrow
  \mbox{Ext}^{1}(L(\overrightarrow{x}_{i}),E)
\longrightarrow \mbox{Ext}^{1}(E,E) \longrightarrow
\mbox{Ext}^{1}(L(\overrightarrow{\omega}),E) \longrightarrow  0 .
 \end{aligned}$$

  \vspace{2mm}

\noindent It follows that
$$\mbox{Hom}(E,E)=\mbox{Hom}(L(\overrightarrow{\omega}),E) = k \
\text{and}\ \mbox{Ext}^{1}(E,E)=0.$$ Hence, $E$ is  exceptional in
$\mbox{coh}(\mathbb{X})$.

As a result, we only need to prove $E[1]$ and
$E(\overrightarrow{x}_{j})$ ($j\neq i$) have the same class in
$K_{0}(\mathbb{X})$ since they are both exceptional in
$\mbox{coh}\mathbb{X}$.

For any $j\notin I \subseteq \{1, 2, 3, 4\}$, we have the following
two exact sequences
$$0 \longrightarrow L(\overrightarrow{\omega}) \longrightarrow L(\overrightarrow{\omega}+\overrightarrow{x}_{j}) \longrightarrow S_{j}\longrightarrow 0 $$
and
$$0 \longrightarrow L(\overrightarrow{\omega}+\sum\limits_{i\in I}\overrightarrow{x}_{i})
\longrightarrow
L(\overrightarrow{\omega}+\overrightarrow{x}_{j}+\sum\limits_{i\in
I}\overrightarrow{x}_{i} ) \longrightarrow S_{j}(\sum\limits_{i\in
I}\overrightarrow{x}_{i} )\longrightarrow 0,$$ where $S_{j}$ denotes
the unique simple sheaf concentrated in the point $x_{j}$ with
$\mbox{Hom}(L ,S_{j})= k$. By noticing that $S_{j}(\sum\limits_{i\in
I}\overrightarrow{x}_{i})=S_{j}$ since $j\notin I$, we obtain
$$[L(\overrightarrow{\omega}+\sum\limits_{i\in
I}\overrightarrow{x}_{i})]+[L(\overrightarrow{\omega}+\overrightarrow{x}_{j})]=[L(\overrightarrow{\omega}+\overrightarrow{x}_{j}+\sum\limits_{i\in
I}\overrightarrow{x}_{i} )]+[L(\overrightarrow{\omega})].$$ Then
from the exact sequence
$$ \xymatrix@C=0.3cm{
  0 \ar[rr]^{} && E \ar[rr]^{\alpha} && IE \ar[rr]^{ \beta} && E[1] \ar[rr]^{ } && 0 },$$
where $IE=L(\overrightarrow{x}_{i})\bigoplus(\bigoplus\limits_{j\neq
i}L(\overrightarrow{\omega}+\overrightarrow{x}_{j}))$, we obtain
that $[E[1]]=[IE]-[E] =\sum\limits_{j\neq
i}[L(\overrightarrow{\omega}+\overrightarrow{x}_{j})]-[L(\omega)]
=[L(\overrightarrow{\omega}+\sum\limits_{j\neq
i}\overrightarrow{x}_{j})]+[L(\omega)]=[L(\overrightarrow{x}_{i}+\overrightarrow{c})]+[L(\overrightarrow{\omega})]$,
and $[E(\overrightarrow{x}_{j})]
=[L(\omega+\overrightarrow{x}_{j})]+[L(\overrightarrow{x}_{i}+\overrightarrow{x}_{j})].$
Now considering the following two exact sequences in
$\mbox{coh}\mathbb{X}$:
$$0 \longrightarrow L(\overrightarrow{\omega})
\longrightarrow L(\overrightarrow{\omega}+\overrightarrow{x}_{j}) \longrightarrow S_{j}\longrightarrow 0 $$
and
$$0 \longrightarrow L(\overrightarrow{x}_{i}+\overrightarrow{x}_{j})
\longrightarrow L(\overrightarrow{x}_{i}+\overrightarrow{c})  \longrightarrow S'\longrightarrow 0, $$
we obtain
$S'=S_{j}(\overrightarrow{x}_{i}+\overrightarrow{x}_{j}-\overrightarrow{\omega})=
S_{j}$ for each $j\neq i$, which implies
$[E[1]]=[E(\overrightarrow{x}_{j})]$. $\hfill\square$

{\bf Remark} \ The proof of the first statement in Proposition 5.2 is an instance of mutations for an exceptional pair, compare [2].

{\bf Corollary 5.3} For each line bundle $L,$  we have

$$E_{L}\langle \overrightarrow{x}_{i} \rangle[n]
=E_{L}\langle\overrightarrow{x}_{i} \rangle(\sum\limits_{j=1, j\neq
i}^4 k_{j} \overrightarrow{x}_{j}),$$ \noindent  where $k_{j}\in
\mathbb{Z} $ satisfying $\sum\limits_{j=1, j\neq i}^4{k_{j}}=n$.
Then $\mu(E_{L}\langle \overrightarrow{x}_{i}
\rangle[n])=n+\mu(E_{L}\langle\overrightarrow{x}_{i} \rangle).$

$\bf{Proof}$: We prove the result by induction.

For $n=1$, there are two possibilities for $\sum\limits_{j=1, j\neq
i}^4{k_{j}}=1$. In the first case, there exists $j\in
\{1,2,3,4\}\backslash \{i\}$ such that $\sum\limits_{j=1, j\neq i}^4
k_{j} \overrightarrow{x}_{j}=\overrightarrow{x}_{j},$ we have
$[E_{L}\langle \overrightarrow{x}_{i} \rangle[1]]
=[L(\overrightarrow{x}_{i}+\overrightarrow{c})]+[L(\overrightarrow{\omega})]
=[E_{L}\langle\overrightarrow{x}_{i}
\rangle(\overrightarrow{x}_{j})]$ from Proposition 5.2. In the other
case, $\sum\limits_{j=1, j\neq i}^4 k_{j}
\overrightarrow{x}_{j}=\sum\limits_{j\neq
i}\overrightarrow{x}_{j}-\overrightarrow{c},$ we have $[E_{L}\langle
\overrightarrow{x}_{i} \rangle(\sum\limits_{j\neq
i}\overrightarrow{x}_{j}-\overrightarrow{c} )]
=[L(\overrightarrow{\omega}+\sum\limits_{j\neq
i}\overrightarrow{x}_{j}-\overrightarrow{c})]
+[L(\overrightarrow{x}_{i}+\sum\limits_{j\neq
i}\overrightarrow{x}_{j}-\overrightarrow{c} )]
=[L(\overrightarrow{x}_{i})]+[L(\overrightarrow{\omega}+\overrightarrow{c})].$
Now considering the following two exact sequences in
$\mbox{coh}\mathbb{X}$:
$$0 \longrightarrow L(\overrightarrow{\omega}) \longrightarrow L(\overrightarrow{\omega}+\overrightarrow{c})
\longrightarrow S(\overrightarrow{\omega})\longrightarrow 0 $$
and
$$0 \longrightarrow L(\overrightarrow{x}_{i}) \longrightarrow L(\overrightarrow{x}_{i}+\overrightarrow{c})
\longrightarrow S(\overrightarrow{x}_{i})\longrightarrow 0, $$
we obtain
$[S(\overrightarrow{\omega})]=[S(\overrightarrow{x}_{i})]=[S]$,
which implies $[E_{L}\langle \overrightarrow{x}_{i}
\rangle(\sum\limits_{j\neq
i}\overrightarrow{x}_{j}-\overrightarrow{c} )] =[E_{L}\langle
\overrightarrow{x}_{i} \rangle[1]]$. Hence, $E_{L}\langle
\overrightarrow{x}_{i} \rangle[1]
=E_{L}\langle\overrightarrow{x}_{i}
\rangle(\sum\limits_{\sum{k_{j}}=1, j\neq i} k_{j}
\overrightarrow{x}_{j}). $

Suppose for $n=k$, the result holds. For $n=k+1,$ noticing that
$E_{L}\langle \overrightarrow{x}_{i}\rangle(\overrightarrow{x})=
E_{L(\overrightarrow{x})}\langle \overrightarrow{x}_{i}\rangle$, we
get $$E_{L}\langle \overrightarrow{x}_{i} \rangle[k+1]=E_{L}\langle
\overrightarrow{x}_{i} \rangle[k][1]=
E_{L}\langle\overrightarrow{x}_{i} \rangle(\sum\limits_{j=1, j\neq
i}^4 k_{j} \overrightarrow{x}_{j}),$$ where $k_{j}\in \mathbb{Z} $
satisfying $\sum\limits_{j=1, j\neq i}^4{k_{j}}=k+1$. This finishes
the proof. $\hfill\square$

{\bf Lemma 5.4} For any indecomposable objects $X,Y\in
\mbox{vect}\mathbb{X}$ , we have two exact sequences:
 $$(1)\ 0 \to \mbox{Hom}(X, Y[-1]) \to \mbox{Hom}(X, PY ) \to
\mbox{Hom}(X, Y ) \to \underline{\mbox{Hom}}(X, Y ) \to 0 ;$$
$$(2)\ \ \  0\to \mbox{Hom}(X[1], Y)\to
\mbox{Hom}(IX, Y )\to\mbox{Hom}(X, Y )\to \underline{\mbox{Hom}}(X,
Y )\to  0\ \  .$$

$\bf{Proof}$: Applying $\mbox{Hom}(X,-)$ to the distinguished exact
sequence
$$ \xymatrix@C=0.3cm{
  0 \ar[rr]^{} && Y[-1] \ar[rr]^{\alpha} && PY \ar[rr]^{ \beta} && Y \ar[rr]^{ } && 0 },$$
we obtain an exact sequence: $$\xymatrix@C=0.3cm{
  0 \ar[rr]^{} &&\mbox{Hom}(X, Y[-1]) \ar[rr]^{\alpha_{\ast}} && \mbox{Hom}(X, PY ) \ar[rr]^{\beta_{\ast}} && \mbox{Hom}(X, Y )
  }.$$
For each $\varphi \in  \mbox{Hom}(X, Y ),$ we have $
\overline{\varphi}=0 \in  \underline{\mbox{Hom}}(X, Y )$  if and
only if $\varphi$ factors through a direct sum of line bundles and
this is equivalent to the fact that  $\varphi$ factors through $PY$
since $\beta: PY\longrightarrow Y$ is distinguished surjective. So
$\mbox{Cok}(\beta_{\ast})=\mbox{Hom}(X, Y )/\mbox{Im}\beta_{\ast}=
\underline{\mbox{Hom}}(X, Y )$. We obtain the exact sequence (1).

Similarly, applying $\mbox{Hom}(-, Y)$ to the distinguished exact
sequence
$$ \xymatrix@C=0.3cm{
  0 \ar[rr]^{} && X \ar[rr]^{\alpha} && IX \ar[rr]^{ \beta} && X[1] \ar[rr]^{} && 0 },$$
we obtain the exact sequence (2).
$\hfill\square$\\

Now we describe all the exceptional objects in
$\underline{\mbox{vect}}\mathbb{X}$.

{\bf Theorem 5.5} A vector bundle $E$ is exceptional in
$\underline{\mbox{vect}}\mathbb{X}$ if and only if  $E$ is an
Auslander bundles or a vector bundle with $\mbox{rk}(E)=p$ and
$\mu(E)=\frac{q}{p}\notin \mathbb{Z},(p,q)=1$.

$\bf{Proof}$: Let $E$ be an indecomposable object in
$\underline{\mbox{vect}}\mathbb{X}$. Then for any $n\geq 2$ or
$n\leq -1,$ we have $\underline{\mbox{Hom}}(E,
E[n])=D\underline{\mbox{Hom}}( E[n-1],
E(\overrightarrow{\omega}))=0$. So $E$ is exceptional in
$\underline{\mbox{vect}}\mathbb{X}$ if and only if
$\underline{\mbox{Hom}}( E, E(\overrightarrow{\omega}))=0$ and
$\underline{\mbox{Hom}}(E, E)= k$.

Notice that if $E(\overrightarrow{\omega})=E$, then
$\underline{\mbox{Hom}}( E, E(\overrightarrow{\omega}))\neq 0$, and
then $E$ can not be exceptional. We only need to consider the
following two cases:

Case 1: $\mu(E)=\frac{q}{p}\notin \mathbb{Z},(p,q)=1$. Suppose
$\mbox{deg}(E)=qr$ and $\mbox{rk}(E)=pr$, then $\mbox{Hom}(IE,
E(\overrightarrow{\omega}))=0$.
 According to the knowledge of the hammocks, we have
$$\mbox{dim }\underline{\mbox{Hom}}( E,
E(\overrightarrow{\omega}))=\mbox{dim Hom}(E,
E(\overrightarrow{\omega}))= [\frac{r}{2}],$$ \noindent here,
$[\frac{r}{2}]$ means  the integral part of $\frac{r}{2}$. So
$\mbox{dim }\underline{\mbox{Hom}}( E,
E(\overrightarrow{\omega}))=0$ if and only if $ r=1$. In this case,
$\mbox{dim } \underline{\mbox{Hom}}(E, E)=\mbox{dim Hom}(E,
E)=[\frac{r+1}{2}]=1$. Hence $E$ is exceptional. So if
$\mu(E)=\frac{q}{p}\notin \mathbb{Z},(p,q)=1$, then $E$ is
exceptional if and only if  $\mbox{rk}E =p$.

Case 2: $\mu(E)=n$. Suppose $\mbox{rk}(E)=r$,  we have $\mbox{dim
Hom}(E, E)=[\frac{r+1}{2}] $ and $\mbox{dim Hom}(IE, E)=
\frac{1+(-1)^{r+1}}{2}  $, so $\mbox{dim }\underline{\mbox{Hom}}(E,
E)=\mbox{dim Hom}(E, E)-\mbox{dim Hom}(IE,
E)=[\frac{r+1}{2}]-\frac{1+(-1)^{r+1}}{2}$. Then $\mbox{dim
}\underline{\mbox{Hom}}(E, E)=1 $ if and only if $r=2$. Hence
$\mbox{dim Hom}(E, E(\overrightarrow{\omega}))=[\frac{r}{2}]=1$ and
$\mbox{dim Hom}(IE,
E(\overrightarrow{\omega}))=\frac{1+(-1)^{r}}{2}=1$. Then we get
$\mbox{dim } \underline{\mbox{Hom}}( E, E(\overrightarrow{\omega}))
=\mbox{dim Hom}(E, E(\overrightarrow{\omega}))-\mbox{dim Hom}(IE,
E(\overrightarrow{\omega}))=0$. Hence $E$ is exceptional. Therefore,
if $\mu(E)=n$, then $E$ is exceptional if and only if
$E(\overrightarrow{\omega})\not= E$ and $r=2$, that is, $E$ is an
Auslander bundle. $\hfill\square$

\section{Tilting objects}

\hspace{4mm} Recall that an object $T$ in a triangulated category
$\mathcal{D}$ is called a tilting object if
$\mbox{Hom}_{\mathcal{D}}(T, T[n])=0$ for $n\neq 0$, and $T$
generates $\mathcal{D}$ as a triangulated category, that is, the
smallest subcategory $\langle T \rangle$ of $\mathcal{D}$, closed
under shift [1] and [-1], direct sums and direct summands, third
terms of triangles, equal to $\mathcal{D}$.

As shown in [6],
$T=\bigoplus\limits_{0\leq\overrightarrow{x}\leq\overrightarrow{c}}\mathcal{O}(\overrightarrow{x}+\overrightarrow{x}_{1})$
is a tilting object in $D^{b}(\mbox{coh}\mathbb{X})$. Then under the
auto-equivalence $\rho$ of $D^{b}(\mbox{coh}\mathbb{X})$ which
acting on slopes $q$ by $q\mapsto \frac{q}{1+q}$, the image $\rho T$
is a tilting object in $D^{b}(\mbox{coh}\mathbb{X})$. Then $\rho T$
is also a tilting object in  $\underline{\mbox{vect}}\mathbb{X}$
since all of the indecomposable direct summands of $\rho T$  are not
line bundles.

Our motivation is to look for some tilting objects directly in
$\underline{\mbox{vect}}\mathbb{X}$.

From now on, we are going to construct a tilting object in $
\underline{\mbox{vect}} \mathbb{X} $, and  always use the notation
$E=E_{\mathcal{O}},
E_{i}=E_{\mathcal{O}}\langle\overrightarrow{x}_{i}\rangle$ .

{\bf Lemma 6.1} For any $n \in \mathbb{Z}$, we have

 (1) $\underline{\mbox{Hom}}(E_{i},E[n])=0;$

 (2) $\underline{\mbox{Hom}}(E,E_{i}[n])=\delta_{n,0}k.$

{\bf{Proof:}} (1)  If $n\leq 0$, then $\mu(E[n]) \leq \mu(E) =0<
\frac{1}{2} =\mu(E_{i}) $ implies $ \mbox{Hom} (E_{i},E[n])=0$, so
$\underline{\mbox{Hom}}(E_{i},E[n])=0$. If $n\geq 2$, then
$\mu(E[n-1]) \geq \mu(E[1]) =\frac{4}{3}> \frac{1}{2} =\mu
(E_{i}(\overrightarrow{\omega}))) $, so $
\mbox{Hom}(E[n-1],E_{i}(\overrightarrow{\omega}))=0$. By Serre
duality,
 $\underline{\mbox{Hom}}(E_{i},E[n])=D\underline{\mbox{Hom}} (E[n-1],E_{i}(\overrightarrow{\omega}))=0$.
 If $n=1$, then $\mbox{dim Hom}(E ,E_{i}(\overrightarrow{\omega}))=\mbox{deg}( E_{i}(\overrightarrow{\omega}) )=1$,
and $\mbox{dim Hom}(IE ,E_{i}(\overrightarrow{\omega}))=\mbox{dim
Hom}(\mathcal{O},E_{i}(\overrightarrow{\omega}))=1$, which implies
$\underline{\mbox{Hom}} (E ,E_{i}(\overrightarrow{\omega}))=0$. By
Serre duality, we have
$\underline{\mbox{Hom}}(E_{i},E[1])=D\underline{\mbox{Hom}} (E,
E_{i}(\overrightarrow{\omega}))=0$.

(2)  For any $n\neq 0$, we obtain
$\underline{\mbox{Hom}}(E,E_{i}[n])=0$ as shown in (1). For $n=0$,
we have $\mbox{dim Hom}(E ,E_{i})=\mbox{deg} (E_{i})=1$ and
$\mbox{dim Hom}(IE ,E_{i} )=\mbox{dim Hom}(\mathcal{O} ,E_{i} )=0$,
so  $\underline{\mbox{Hom}}(E,E_{i} )= 1.$ $\hfill\square$

 Applying  the functor $\mbox{Hom}(-,\mathcal{O}(\overrightarrow{\omega}))$ to the Auslander-Reiten sequence
$$ \xi:
0 \longrightarrow \mathcal{O}(\overrightarrow{x}_{j})\longrightarrow
E(\overrightarrow{\omega} + \overrightarrow{x}_{j}) \longrightarrow
\mathcal{O}(\overrightarrow{\omega}+\overrightarrow{x}_{j})\longrightarrow
0,$$ we obtain an exact sequence
$$0 \longrightarrow \mbox{Ext}^{1}(\mathcal{O}(\overrightarrow{\omega} +
\overrightarrow{x}_{j}),\mathcal{O}( \overrightarrow{\omega} ))\longrightarrow \mbox{Ext}^{1}(E(\overrightarrow{\omega} +
\overrightarrow{x}_{j}),\mathcal{O}( \overrightarrow{\omega} )) $$

\hspace{50mm}$\longrightarrow \mbox{Ext}^{1}(\mathcal{O}(
\overrightarrow{x}_{j}),\mathcal{O}( \overrightarrow{\omega}
))\longrightarrow 0.$

\noindent Then $\mbox{Ext}^{1}(\mathcal{O}(\overrightarrow{\omega} +
\overrightarrow{x}_{j}),\mathcal{O}(\overrightarrow{ \omega }))=
S_{\overrightarrow{\omega }+ \overrightarrow{x}_{j}}=0$ and
$\mbox{Ext}^{1}(\mathcal{O}( \overrightarrow{x}_{j}),\mathcal{O}(
\overrightarrow{\omega} )) = S_{  \overrightarrow{x}_{j}}= k$, which
imply $\mbox{Ext}^{1}(E(\overrightarrow{\omega} +
\overrightarrow{x}_{j}),\mathcal{O}(\overrightarrow{ \omega })) =
k$. Hence, there is a vector bundle $F$ fitting into the following
exact sequence
$$\zeta:   0 \longrightarrow \mathcal{O}(\overrightarrow{\omega})\longrightarrow F
\longrightarrow E(\overrightarrow{\omega} + \overrightarrow{x}_{j})\longrightarrow 0.$$
It is easy to see that  $\mbox{deg}(F)=2$, and $\mbox{rk}(F)=3$.
Moreover, $F$ is indecomposable since there is no line bundle $L$
satisfying  $\mbox{Hom}(\mathcal{O}( \overrightarrow{\omega} ),
L)\neq 0$ and $\mbox{Hom}(L, E(\overrightarrow{\omega} +
\overrightarrow{x_{j}}))\neq 0 $.

{\bf Theorem 6.2}  Let $T=E\bigoplus(\bigoplus\limits_{
i=1}^{4}E_{i} )\bigoplus F$, then

(1) $T$ is a tilting object in $\underline{\mbox{vect}}\mathbb{X}$;

(2) $\mbox{End}(T)$ is a canonical algebra of type (2,2,2,2).

{\bf{Proof: }} (1) First we show that
$\underline{\mbox{Hom}}(T,T[n])=0$ for any $n\neq 0$ and  that $E,
E_{1}, E_{2}, E_{3}$, $E_{4}$ ,$F$ forms an exceptional sequence.

By comparing  slopes and using the Serre duality, we obtain:

(i) \ For  $n\neq 0$,
$\underline{\mbox{Hom}}(E_{i},F[n])=D\underline{\mbox{Hom}}
(F[n-1],E_{i}(\overrightarrow{\omega}))=0$ and
$\underline{\mbox{Hom}}(E ,F[n])=D\underline{\mbox{Hom}} (F[n-1],E
(\overrightarrow{\omega}))=0$.

(ii) \ For  $n\neq 1$, $\underline{\mbox{Hom}}(F,
E_{i}[n])=D\underline{\mbox{Hom}}
(E_{i}[n-1],F(\overrightarrow{\omega}))=0$ and
$\underline{\mbox{Hom}}(F ,E[n])=D\underline{\mbox{Hom}}
(E[n-1],F(\overrightarrow{\omega}))=0$. According to Lemma 6.1, in
order to prove $T$ is a tilting object in
$\underline{\mbox{vect}}\mathbb{X}$, we only need to show that $
\underline{\mbox{Hom}}(F, E_{i}[1])=0$  for each $i=1,2,3,4$ and
$\underline{\mbox{Hom}}(F, E[1])=0$.

Applying $\mbox{Hom}(E_{i},-)$ to the exact sequence
$$\zeta(\overrightarrow{\omega}):   0 \longrightarrow \mathcal{O}
\longrightarrow F(\overrightarrow{\omega}) \longrightarrow E(\overrightarrow{x}_{j})\longrightarrow 0,$$
we obtain an exact sequence:
$$0\longrightarrow \mbox{Hom}(E_{i},F(\overrightarrow{\omega}))
\longrightarrow
\mbox{Hom}(E_{i},E(\overrightarrow{x}_{j}))\longrightarrow
\mbox{Ext}^{1}(E_{i},\mathcal{O})\longrightarrow 0.$$ Then
$\mbox{Hom}(E_{i}, F(\overrightarrow{\omega}))=0$ provided by
$$\mbox{dim Ext}^{1}(E_{i}, \mathcal{O})= \mbox{dim
DHom}(\mathcal{O}(\overrightarrow{\omega}), E_{i})=1$$ and
$$\mbox{dim Hom}(E_{i}, E(\overrightarrow{x}_{j}))=\mbox{dim
Hom}(E_{i}(-\overrightarrow{x}_{j}),
E)=-\mbox{deg}(E_{i}(-\overrightarrow{x}_{j}))=1.$$ Therefore,
$\underline{\mbox{Hom}}(F, E_{i}[1])=\mbox{D}\underline{\mbox{Hom}}
(E_{i} ,F(\overrightarrow{\omega}))=0$  for each $i=1,2,3,4$.

Applying $\mbox{Hom}(\mathcal{O},-)$ to the exact sequence $\zeta$
and $\xi$, we obtain two exact sequences:
$$0\longrightarrow \mbox{Hom}(\mathcal{O},F )
\longrightarrow \mbox{Hom}(\mathcal{O},E(\overrightarrow{\omega}+\overrightarrow{x}_{j}))
\longrightarrow
\mbox{Ext}^{1}(\mathcal{O},\mathcal{O}(\overrightarrow{\omega}))\longrightarrow
0,$$ and
$$0\to \mbox{Hom}(\mathcal{O},\mathcal{O}( \overrightarrow{x}_{j})
)\to
\mbox{Hom}(\mathcal{O},E(\overrightarrow{\omega}+\overrightarrow{x}_{j}))\to
\mbox{Hom}(\mathcal{O},\mathcal{O}(\overrightarrow{\omega}+\overrightarrow{x}_{j}))\to
0.$$ So it is easy to show that $\mbox{dim Hom}(\mathcal{O},F )=0$.
Then it follows that $\mbox{dim
Hom}(\mathcal{O}(\overrightarrow{\omega}),F )=2$ since $\mbox{dim
Hom}(\mathcal{O}\bigoplus \mathcal{O}(\overrightarrow{\omega}), F
)=\mbox{deg}(F) =2$. Therefore, $$\mbox{dim
Hom}(IE,F(\overrightarrow{\omega}))=\mbox{dim
Hom}(\mathcal{O},F(\overrightarrow{\omega}))=\mbox{dim
Hom}(\mathcal{O}(\overrightarrow{\omega}),F)=2.$$ By the fact that
$\mbox{dim Hom}(E
,F(\overrightarrow{\omega}))=\mbox{deg}(F(\overrightarrow{\omega}))=2$,
we have $$\underline{\mbox{Hom}}(F,
E[1])=\mbox{D}\underline{\mbox{Hom}}(E,F(\overrightarrow{\omega}) )=
0.$$

Next, we need to show $T$ generates
$\underline{\mbox{vect}}\mathbb{X}$. Since $E, E_{1}, E_{2}, E_{3},
E_{4} , F$ is an exceptional sequence,
 it suffices to prove that for each indecomposable vector bundle $X$ with $\mbox{rk}(X)\geq 2$,
 there exists some $n\in \mathbb{Z}$, such that $\mbox{\underline{Hom}}(T, X[n])\neq 0$.

Indeed, if we fix an indecomposable vector bundle $X$ with
$\mbox{rk}(X)\geq 2$, then there exists some $m\in \mathbb{Z}$ such
that  $m\leq\mu(X)<m+1$. By Corollary 4.7,  we have
$-\frac{m}{2m+1}=\mu(E_{\mathcal
{O}(m\overrightarrow{x}_{1})})-m-\frac{m}{2m+1}=\mu(E_{\mathcal
{O}(m\overrightarrow{x}_{1})}[-m])\leq\mu(X[-m])<\mu(E_{\mathcal
{O}((m+1)\overrightarrow{x}_{1})}[-m])=\mu(E_{\mathcal
{O}((m+1)\overrightarrow{x}_{1})})-m-\frac{m}{2m+1}=\frac{m+1}{2m+1}.$
Moreover, if $-\frac{m}{2m+1}\leq\mu(X[-m])< 0$, then by Proposition
4.6, we have $\frac{m+1}{2m+1}\leq\mu(X[-m][1])< \frac{4}{3}$.
Hence, there exists a suitable integer $n_{1}\in \mathbb{Z}$ such
that $0\leq\mu(X[n_{1}])<\frac{4}{3}.$ So we only need to show that
for each indecomposable vector bundle $X$ with $\mbox{rk}(X)\geq 2$
and $0\leq\mu(X)<\frac{4}{3}$, there exists some $n\in \mathbb{Z}$,
such that $\mbox{\underline{Hom}}(T, X[n])\neq 0$.

(i) If $0< \mu(X)< 1$, then
$$\underline{\mbox{Hom}}(E,X)=\mbox{Hom}(\mathcal{O}(\overrightarrow{\omega})\bigoplus
\mathcal{O},X)-
\mbox{Hom}(IE,X)=\mbox{Hom}(\mathcal{O}(\overrightarrow{\omega}),X).$$
If $\mbox{Hom}(\mathcal{O}(\overrightarrow{\omega}),X)\neq 0$, then
$\underline{\mbox{Hom}}(E,X)\neq 0.$ If
$\mbox{Hom}(\mathcal{O}(\overrightarrow{\omega}),X)=0$, then
$\mbox{Ext}^{1}(X,\mathcal{O})=0.$ Applying $\mbox{Hom}(X,-)$ to the
exact sequence
$$0\longrightarrow \mathcal{O}\longrightarrow F(\overrightarrow{\omega})\longrightarrow E(\overrightarrow{x}_{j})\longrightarrow 0,$$
we have
$$0\longrightarrow \mbox{Hom}(X,F(\overrightarrow{\omega}))
\longrightarrow \mbox{Hom}(X,E(\overrightarrow{x}_{j}))\longrightarrow \mbox{Ext}^{1}(X,\mathcal{O})=0.$$
So $\mbox{dim }\mbox{Hom}(X,F(\overrightarrow{\omega}))= \mbox{dim
}\mbox{Hom}(X,E(\overrightarrow{x}_{j}))=-\mbox{deg}(X(-\overrightarrow{x}_{j}))=\mbox{rk}(X)-\mbox{deg}(X)=\mbox{rk}(X)(1-\mu(X))>
0.$ Hence,
$\underline{\mbox{Hom}}(F,X[1])=\mbox{D}\underline{\mbox{Hom}}(X,F(\overrightarrow{\omega}))=\mbox{D}\mbox{Hom}(X,F(\overrightarrow{\omega}))\neq
0.$

(ii) If $\mu(X)=1$, applying $\mbox{Hom}(-,X)$ to the exact sequence
$$0\longrightarrow \mathcal{O}(\overrightarrow{\omega})
\longrightarrow F\longrightarrow
E(\overrightarrow{\omega}+\overrightarrow{x}_{j})\longrightarrow
0,$$ we have
$$\begin{aligned}
0&\longrightarrow
\mbox{Hom}(E(\overrightarrow{\omega}+\overrightarrow{x}_{j}),X)\longrightarrow
\mbox{Hom}(F,X)\longrightarrow
\mbox{Hom}(\mathcal{O}(\overrightarrow{\omega}),X)
\\&\longrightarrow\mbox{Ext}^{1}(E(\overrightarrow{\omega}+\overrightarrow{x}_{j}),X)\longrightarrow 0.
\end{aligned}$$
Notice that
$$\begin{aligned}
&\mbox{dim
}\mbox{Hom}(E(\overrightarrow{\omega}+\overrightarrow{x}_{j}),X)-\mbox{dim
}\mbox{Ext}^{1}(E(\overrightarrow{\omega}+\overrightarrow{x}_{j}),X)\\&=
\langle
[E(\overrightarrow{\omega}+\overrightarrow{x}_{j})],[X]\rangle\\&=\langle
[\mathcal{O}(\overrightarrow{\omega}+\overrightarrow{x}_{j})]+[\mathcal{O}(\overrightarrow{x}_{j})],[X]\rangle
\\&=\mbox{rk}(\mathcal{O}(\overrightarrow{x}_{j}))\mbox{deg}(X)-\mbox{deg}(\mathcal{O}(\overrightarrow{x}_{j}))\mbox{rk}(X)\\&=\mbox{deg}(X)-\mbox{rk}(X)=0.
\end{aligned}$$
We get $\mbox{dim }\mbox{Hom}(F,X)=\mbox{dim
}\mbox{Hom}(\mathcal{O}(\overrightarrow{\omega}),X)\geq 1.$
Moreover, since
$$IF=\mathcal{O}(\overrightarrow{c})\bigoplus(\bigoplus\limits_{i=1}^{4}\mathcal{O}(\overrightarrow{\omega}+\overrightarrow{x}_{i})),$$
we get $$\mbox{dim }\mbox{Hom}(IF,X)=\mbox{dim
}\mbox{Hom}(\bigoplus\limits_{i=1}^{4}\mathcal{O}(\overrightarrow{\omega}+\overrightarrow{x}_{i}),X)\leq
1.$$ Hence, $ \underline{\mbox{Hom}}(F, X)=0$ if and only if
$X=E_{\mathcal{O}(\overrightarrow{x}_{i})}.$ In this case,
$$\mbox{dim
Hom}(E_{i},E_{\mathcal{O}(\overrightarrow{x}_{i})})=\mbox{dim
Hom}(E_{i}(-\overrightarrow{x}_{i}),E)=-\mbox{deg}E_{i}(-\overrightarrow{x}_{i})=1$$
and $$\mbox{Hom}(IE_{i},E_{\mathcal{O}(\overrightarrow{x}_{i})})=0$$
imply that
$\underline{\mbox{Hom}}(E_{i},E_{\mathcal{O}(\overrightarrow{x}_{i})})\neq
0.$

(iii) If $\mu(X)=0$, then by Serre duality,
$$\begin{aligned}
&\ \ \ \ \mbox{dim }\underline{\mbox{Hom}}(E_{i},X[1])=\mbox{dim
}\mbox{D}\underline{\mbox{Hom}}(X,E_{i}(\overrightarrow{\omega}))\\&=\mbox{dim
}\mbox{DHom}(X,E_{i}(\overrightarrow{\omega}))-\mbox{dim
}\mbox{DHom}(X,PE_{i}(\overrightarrow{\omega})).
\end{aligned}$$ Applying
$\mbox{Hom}(X,-)$ to the exact sequence
$$0\longrightarrow \mathcal{O}\longrightarrow E_{i}(\overrightarrow{\omega})
\longrightarrow \mathcal{O}(\overrightarrow{\omega}+\overrightarrow{x}_{i})\longrightarrow 0,$$
we obtain
$$\begin{aligned}
0&\longrightarrow \mbox{Hom}(X,\mathcal{O})\longrightarrow
\mbox{Hom}(X,E_{i}(\overrightarrow{\omega})) \longrightarrow
\mbox{Hom}(X,\mathcal{O}(\overrightarrow{\omega}+\overrightarrow{x}_{i}))\\&
\longrightarrow\mbox{Ext}^{1}(X,\mathcal{O})\longrightarrow 0.
\end{aligned}$$
Notice that $P(E_{i}(\overrightarrow{\omega}))=\mathcal{O}\bigoplus
(\bigoplus\limits_{j\neq
i}\mathcal{O}(\overrightarrow{\omega}+\overrightarrow{x}_{i}-\overrightarrow{x}_{j}))$.
We obtain that
$$\begin{aligned}
&\ \ \ \ \mbox{dim
}\underline{\mbox{Hom}}(X,E_{i}(\overrightarrow{\omega}))\\
&=\mbox{dim }\mbox{Hom}(X,E_{i}(\overrightarrow{\omega}))-\mbox{dim }\mbox{Hom}(X,P(E_{i}(\overrightarrow{\omega})))\\
&=\mbox{dim
}\mbox{Hom}(X,\mathcal{O}({\overrightarrow{\omega}}+\overrightarrow{x}_{i}))-\sum\limits_{j\neq
i}\mbox{dim }\mbox{Hom}(X,\mathcal{O}(\overrightarrow{\omega}+\overrightarrow{x}_{i}-\overrightarrow{x}_{j}))\\
&\ \ \ \ -\mbox{dim }\mbox{Ext}^{1}(X,\mathcal{O})\geq
0.\end{aligned}$$ Moreover,
$\underline{\mbox{Hom}}(X,E_{i}(\overrightarrow{\omega})=0$ if and
only if $X=E$ or
$X=E_{\mathcal{O}(\overrightarrow{\omega}+\overrightarrow{x}_{i}-\overrightarrow{x}_{j})}$
for some $j\neq i$. So we only need to consider the case
$X=E_{\mathcal{O}(\overrightarrow{\omega}+\overrightarrow{x}_{i}-\overrightarrow{x}_{j})}$.
Since
$PF=\mathcal{O}(\overrightarrow{\omega})^{2}\bigoplus(\bigoplus\limits_{i\neq
j}\mathcal{O}(\overrightarrow{x}_{i}-\overrightarrow{x}_{j}))$ and
$\mbox{deg}(F(\overrightarrow{\omega}))=2,$ we have $\mbox{dim
}\underline{\mbox{Hom}}(F,E_{\mathcal{O}(\overrightarrow{\omega}+\overrightarrow{x}_{i}-\overrightarrow{x}_{j})}[1])=
\mbox{dim
}\mbox{D}\underline{\mbox{Hom}}(E_{\mathcal{O}(\overrightarrow{\omega}+\overrightarrow{x}_{i}-\overrightarrow{x}_{j})},F(\overrightarrow{\omega}))=1.$

(iv) If $1<\mu(X)<\frac{4}{3}$, then $\mbox{Hom}(X,
\mathcal{O}(\overrightarrow{c}))\neq 0.$ Notice that
$$\underline{\mbox{Hom}}(F,X)=\mbox{D}\underline{\mbox{Hom}}(X(\overrightarrow{\omega}),F[1])=
\mbox{D}\underline{\mbox{Hom}}(X(\overrightarrow{\omega}),E_{\mathcal{O}(\overrightarrow{c}+\overrightarrow{\omega})}).$$
It follows that
$$\begin{aligned}
&\ \ \ \ \mbox{dim }\underline{\mbox{Hom}}(F,X)\\&= \mbox{dim
}\mbox{Hom}(X(\overrightarrow{\omega}),E_{\mathcal{O}(\overrightarrow{c}+\overrightarrow{\omega})})-
\mbox{dim
}\mbox{Hom}(X(\overrightarrow{\omega}),PE_{\mathcal{O}(\overrightarrow{c}+\overrightarrow{\omega})})\\&
=\mbox{dim
}\mbox{Hom}(X(\overrightarrow{\omega}),\mathcal{O}(\overrightarrow{c}+\overrightarrow{\omega})\bigoplus
\mathcal{O}(\overrightarrow{c}))-\mbox{dim
}\mbox{Hom}(X(\overrightarrow{\omega}),\mathcal{O}(\overrightarrow{c}))\\&=
\mbox{dim
}\mbox{Hom}(X(\overrightarrow{\omega}),\mathcal{O}(\overrightarrow{c}+\overrightarrow{\omega})\\&=\mbox{dim
}\mbox{Hom}(X,\mathcal{O}(\overrightarrow{c}))\neq 0.
\end{aligned}$$

(2) Since $\mbox{dim Hom}(E,F )=\mbox {deg} (F )=2$ and $$\mbox{dim
Hom}(IE,F )=\mbox{dim Hom}(\mathcal{O},F )=0,$$ we have $\mbox{dim
}\underline{\mbox{Hom}}(E,F )=2$. Furthermore, for each $i=1,2,3,4$,
applying $\mbox{Hom}(-,F)$ to the exact sequence
$$\eta_{i}: 0 \longrightarrow \mathcal{O}(\overrightarrow{\omega})\longrightarrow E_{i}
\longrightarrow \mathcal{O}(\overrightarrow{x}_{i})\longrightarrow
0,$$ we obtain an exact sequence
$$0\longrightarrow \mbox{Hom}( E_{i},F )\longrightarrow \mbox{Hom}(\mathcal{O}( \overrightarrow{\omega}),F)
\longrightarrow \mbox{Ext}^{1}(\mathcal{O}(
\overrightarrow{x}_{i}),F)\longrightarrow 0;$$ and applying
$\mbox{Hom}(\mathcal{O}( \overrightarrow{x}_{i}),-)$ to the exact
sequence $\zeta$, we obtain an exact sequence

\hspace{15mm}$0\longrightarrow \mbox{Hom}(\mathcal{O}(
\overrightarrow{x}_{i}),E(
\overrightarrow{\omega}+\overrightarrow{x}_{j}))\longrightarrow
\mbox{Ext}^{1}(\mathcal{O}(
\overrightarrow{x}_{i}),\mathcal{O}(\overrightarrow{\omega}))$
$$\longrightarrow \mbox{Ext}^{1}(\mathcal{O}( \overrightarrow{x}_{i}),F)\longrightarrow
\mbox{Ext}^{1}(\mathcal{O}( \overrightarrow{x}_{i}),E(
\overrightarrow{x}_{j}+\overrightarrow{\omega}))\longrightarrow 0.$$
Using that $\mbox{Ext}^{1}(\mathcal{O}( \overrightarrow{x}_{i}),E(
\overrightarrow{\omega}+\overrightarrow{x}_{j}))= \mbox{DHom}(E(
\overrightarrow{x}_{j}),\mathcal{O}( \overrightarrow{x}_{i})),$ we
can easily get that $$\mbox{dim Hom}(\mathcal{O}(
\overrightarrow{x}_{i}),E(
\overrightarrow{\omega}+\overrightarrow{x}_{j}))=\mbox{dim
Ext}^{1}(\mathcal{O}( \overrightarrow{x}_{i}),E(
\overrightarrow{\omega}+\overrightarrow{x}_{j}))=\delta_{i,j},$$
which implies that $$\mbox{dim Ext}^{1}(\mathcal{O}(
\overrightarrow{x}_{i}),F)=\mbox{dim Ext}^{1}(\mathcal{O}(
\overrightarrow{x}_{i}),\mathcal{O}(\overrightarrow{\omega}))=\mbox{dim
} S_{ \overrightarrow{x}_{i}}=1.$$ Hence,  we can obtain $\mbox{dim
Hom}( E_{i},F )=1$ since $\mbox{dim
Hom}(\mathcal{O}(\overrightarrow{\omega}),F)=2$. Moreover,
$\mbox{Hom}( IE_{i},F )=0$ implies $\mbox{dim
}\underline{\mbox{Hom}}( E_{i},F )=1$.

Now we describe generators and the relations. Assume  $1\leq i \leq 4$ throughout  the rest of the proof.

Applying  $\mbox{Hom}( - , \mathcal{O}(
\overrightarrow{c})$ to the exact sequence
$$
\xymatrix@C=0.3cm{
  0 \ar[rr]^{} && \mathcal{O}(\overrightarrow{x}_{i}) \ar[rr]^{\eta_{i}\ \ \ \ } &&
  E(\overrightarrow{\omega}+\overrightarrow{x}_{i})\ar[rr] && \mathcal{O}(\overrightarrow{\omega}+\overrightarrow{x}_{i}) \ar[rr]^{ } && 0
  },$$
then  $\eta_{i}$ induces an isomorphism $\mbox{Hom}(  E(\overrightarrow{\omega}+\overrightarrow{x}_{i}), \mathcal{O}( \overrightarrow{c} )\cong \mbox{Hom}(
\mathcal{O}( \overrightarrow{x}_{i}),
\mathcal{O}(\overrightarrow{c}) )$
sending the generator $\theta_{i}$ to  $x_{i}$, that is, the following diagram (I) commutes:

$$
\xymatrix@C=0.3cm{
  \mathcal{O}(\overrightarrow{x}_{i})\ar[rd]_{x_{i}} \ar[rr]^{\eta_{i}} && E(\overrightarrow{\omega}+
  \overrightarrow{x}_{i})\ar[ld]^{\theta_{i}}&&&&F\ar[rd]_{\gamma} \ar[rr]^{\pi_{i} \ \ } && E(\overrightarrow{\omega}+\overrightarrow{x}_{i})\ar[ld]^{\theta_{i}}\\
 &  \mathcal{O}(\overrightarrow{c}) &&&&&&\mathcal{O}(\overrightarrow{c}) &\\
 &(I)&&&&&&(II)
  },$$
Analogously,  applying
 $\mbox{Hom}( - , \mathcal{O}(
\overrightarrow{c})$ to the exact sequence
$$
\xymatrix@C=0.3cm{
  0 \ar[rr]^{} && \mathcal{O}(\overrightarrow{\omega}) \ar[rr] && F \ar[rr]^{\pi_{i}\ \ \ \ \ \ \  }&&
  E(\overrightarrow{\omega}+\overrightarrow{x}_{i})\ar[rr]^{ }   && 0
  },$$
   then  $\pi_{i}$ induces an isomorphism $\mbox{Hom}(  E(\overrightarrow{\omega}+\overrightarrow{x}_{i}),
   \mathcal{O}( \overrightarrow{c} ))\cong \mbox{Hom}(F, \mathcal{O}(\overrightarrow{c}) )$
sending the generator $\theta_{i}$ to  $\gamma$, that is, the above   diagram (II) commutes.

Now from the following commutative diagram,

$$
\xymatrix@C=0.3cm{
   0 \ar[rr]^{} && \mathcal{O}(\overrightarrow{\omega})
   \ar[rr] \ar@{=}[d]&& E \ar[rr]^{\pi  }\ar[d]^{f_{i}}&& \mathcal{O}\ar[rr]^{ }\ar[d]^{x_{i}}   && 0\\
   0 \ar[rr]^{} && \mathcal{O}(\overrightarrow{\omega})
   \ar[rr] \ar@{=}[d]&& E\langle\overrightarrow{x}_{i}\rangle \ar[rr]^{ }\ar[d]^{g_{i}}&&
   \mathcal{O}(\overrightarrow{x}_{i})\ar[d]^{\eta_{i}}\ar[rr]^{ }  \ar[rrdd]^{x_{i}} && 0,\\
    0 \ar[rr]^{} && \mathcal{O}(\overrightarrow{\omega})
    \ar[rr] && F \ar[rr]^{\pi_{i}\ \ \ \ \ \ \  }\ar[rrrrd]_{\gamma}&& E(\overrightarrow{\omega}+\overrightarrow{x}_{i})\ar[rr]^{ }   \ar[rrd]^{\theta_{i}} && 0\\
   &&&&&&&& \mathcal{O}(\overrightarrow{c})
  }$$
where $f_{i}$ and $g_{i}$ are obtained from pullbacks, we have
$$\gamma (g_{i} f_{i})=\theta_{i}\pi_{i} g_{i}
f_{i}=(\theta_{i}\eta_{i}) x_{i}\pi =x^{2}_{i}\pi.$$ Hence, the fact
that $ \{x^{2}_{i}|1\leq i \leq 4 \}$ are pairwise linearly
independent implies that $ \{g_{i} f_{i}|1\leq i \leq 4 \}$ are
pairwise linearly independent since $\pi$ is surjective.

Next, we claim that $\theta_{i}:
E(\overrightarrow{\omega}+\overrightarrow{x}_{i})\longrightarrow
\mathcal{O}( \overrightarrow{c} )$ is surjective. Otherwise,
$\mbox{Im}(\theta_{i})\subseteq \mathcal{O}( \overrightarrow{c} )$
is a line bundle satisfying $\mbox{Hom}(
E(\overrightarrow{\omega}+\overrightarrow{x}_{i}),
\mbox{Im}(\theta_{i}) )\neq 0$,
 which implies $\mbox{Im}(\theta_{i})=\mathcal{O}(\overrightarrow{\omega}+\overrightarrow{x}_{i}) $. And it follows that
 $\mbox{Ker}(\theta_{i})=\mathcal{O}(\overrightarrow{x}_{i})$, which is a contradiction.

  Moreover, assuming the generators of  $\mbox{Hom}(\mathcal{O}( \overrightarrow{\omega}), F )$ are $ h_{1}, h_{2}$,  one checks easily that
  $\mbox{Ker}(\gamma)=\mathcal{O}(\overrightarrow{\omega})^{2}$, and
  there has the following exact sequence:
$$
\xymatrix@C=0.3cm{
  0 \ar[rr]^{} && \mathcal{O}(\overrightarrow{\omega})^{2} \ar[rr]^{ \  \ \ (h_{1}, h_{2})} && F \ar[rr]^{\gamma\ \ }
  &&\mathcal{O}(\overrightarrow{c} )\ar[rr]^{ }   && 0
  }.$$ Then by applying $\mbox{Hom}(E,  -)$ to this exact sequence one obtains that $\gamma$ induces an monomorphism
  $$\gamma: \mbox{Hom}(E,  F)\longrightarrow \mbox{Hom}(E,   \mathcal{O}(\overrightarrow{c} )); \ \ \  g_{i} f_{i} \longmapsto \gamma (g_{i} f_{i}).$$
 Hence the relations  that  $$x^{2}_{3}=x^{2}_{2}+x^{2}_{1}\ \text{and}\   x^{2}_{4}=x^{2}_{2}+\lambda x^{2}_{1}$$ imply that
$$g_{3} f_{3}=g_{2} f_{2}+g_{1} f_{1}\ \text{and}\ g_{4} f_{4}=g_{2} f_{2}+\lambda g_{1}
f_{1}.$$

Therefore, $\mbox{End}(T)$ is the canonical algebra of type
(2,2,2,2) on generators $ \{g_{i} f_{i}|1\leq i \leq 4 \}$ subject
to the relations $$g_{3}f_{3}=g_{2} f_{2}+g_{1} f_{1}\ \text{and}\
g_{4} f_{4}=g_{2} f_{2}+\lambda g_{1} f_{1}. $$ $\hfill\square$

Similarly, we can obtain that:

{\bf Theorem 6.3}  Let $T'=F[-1](\overrightarrow{\omega}) \bigoplus
E \bigoplus (\bigoplus\limits_{ i=1}^{4}E_{i} )$, then

(1) $T'$ is a tilting object in $\underline{\mbox{vect}}\mathbb{X}$;

(2) $\mbox{End}(T')$ is an  algebra given by the following quiver
with relations:

\setlength{\unitlength}{0.07cm}
\begin{picture}(200,40)
\put(61,18){\vector(1,0){16}} \put(61,16){\vector(1,0){16}}
\put(80,18){\vector(3,1){16}} \put(80,16){\vector(4,-1){16}}
\put(80,14){\vector(2,-1){16}} \put(80,20){\vector(3,2){16}} {\small
\put(67,18){\makebox(6,4)[l]{$\alpha_1$}}
\put(67,12){\makebox(6,4)[l]{$\alpha_2$}}
\put(90,28){\makebox(6,4)[l]{$\beta_1$}}
\put(90,22){\makebox(6,4)[l]{$\beta_2$}}
\put(90,14){\makebox(6,4)[l]{$\beta_3$}}
\put(90,8){\makebox(6,4)[l]{$\beta_4$}}}
\put(110,20){$\beta_{i}\alpha_{1}=\beta_{i}\alpha_{2}$,}
\put(110,10){$i\in \{1,2,3,4 \}$;}
\end{picture}

{\bf Remark:} For any two distinct indecomposable vector bundles
$E_{1}, E_{2}$ of rank two, we have $\mu E_{1}, \mu E_{2}\in
\mathbb{Z}\bigcup (\frac{1}{2}+\mathbb{Z})$ and $\mbox{dim\
}\underline{\mbox{Hom}}( E_{1}, E_{2} )\neq 0$ implies $\mu
E_{1}\leq \mu E_{2}\leq \mu (E_{1}[1])$. So one can easily check
that $\mbox{dim}\ \underline{\mbox{Hom}}( E_{1}, E_{2} )\leq 1$ case
by case.  Hence there doesn't exist any tilting object only
consisting of rank two bundles whose endomorphism algebra is a
canonical algebra. In other words, our result Theorem 6.2 can not be
much more simple .

\vspace {2mm}

{\bf Acknowledgements} The authors would like to thank Helmut
Lenzing for his lectures in Xiamen University in 2011 and his
helpful suggestion. In particular, this work is inspired by a talk
given by Helmut Lenzing during his visit in Xiamen. The authors also
thank Xiaowu Chen for his encouragement.

\begin{description}
 \item{ {\bf\large References}}
 \vspace {3mm}
\item{[1]} \ W. Geigle and H. Lenzing. A class of weighted projective curves arising
in representation theory of finite dimensional algebras. Singularities,
representations of algebras, and Vector bundles, Springer Lect. Notes Math., 1273(1987), 265-297.

\item{[2]} \ A. L. Gorodentsev and A. N. Rudakov. Exceptional vector bundles on projective spaces. Duke Math, J., 54(1)(1987), 115-130.

\item{[3]} \ D. Happel. Triangulated categories in tne representation theory of finite-dimensional algebras, London Mathematical Society
Lecture Note Series 119, Cambridge University Press, Cambridge,
1988.

\item{[4]} \ D. Kussin, H. Lenzing, H. Meltzer. Triangle singularities, ADE-chains and weighted projective lines. arXiv: 1203.5505

\item{[5]} \ D. Kussin, H. Lenzing, H. Meltzer. Nilpotent operators and weighted projective lines. arXiv:1002.3797vl

\item{[6]} \ H. Lenzing. Weighted projective lines and applications.
Representations of Algebras and Related Topics, European Mathematical Society, 153-187. DOI: 10.4171/101-1/5

\item{[7]} \  H. Lenzing and H. Meltzer. Sheaves on a weighted projective line of genus one, and representations of a tubular algebra.
In Representations of algebras, Sixth International Conference,
Ottawa 1992. CMS Conf.Proc. 14, 1993, 313-337.

\item{[8]} \ D. Quillen. Higher algebraic $K-$theory. I. In Algebraic $K-$theory, I: Higher $K-$theories.
Lecture Notes in Mathematics 341, Springer-Verlag, Berlin, 1973, 85-147.

\item{[9]} \  C.M. Ringel. Tame algebras and integral quadratic forms.
Lecture Notes in Math. 1099, Springer, 1984.

\item{[10]} \ J. P. Serre, Faisceaux alg\'{e}briques coh\'{e}rents. Ann. of math. (2), 61(1955), 197-278.

\item{[11]} \ K. Ueda. A Remark on a Theorem of math.AG/0511155. arXiv:math. AG/0604361, 2006.

\end{description}

\end{document}